\documentclass[a4paper,11pt]{amsart}

\usepackage{amsfonts}
\usepackage{amsmath,amssymb}
\usepackage[dvipsnames]{xcolor}
\usepackage{tikz}
\usetikzlibrary{decorations.pathreplacing, decorations.markings, shapes}

\usepackage{mathtools}
\usepackage{latexsym}
\usepackage{amsthm}
\usepackage{mathrsfs}
\usepackage{geometry}
\usepackage{graphics}
\usepackage{calrsfs}
\usepackage{caption}
\usepackage{mathrsfs}
 \usepackage{pgfplots}
 \pgfplotsset{compat=newest}

\usepackage{color, colortbl}

\newcommand{\evalat}[1]{\bigr\rvert_{#1}}
\newcommand{\Id}{\mathrm{Id}}
\newcommand{\set}[2]{\ensuremath{\left\{#1 \,\colon #2\right\}}}
\makeatletter
\def\@map#1#2[#3]{\mbox{$#1 \colon #2 \longrightarrow #3$}}
\def\map#1#2{\@ifnextchar [{\@map{#1}{#2}}{\@map{#1}{#2}[#2]}}
\makeatother

\newtheorem{MainTheorem}{Theorem}

\newtheorem{theorem}{Theorem}[section]
\newtheorem{proposition}[theorem]{Proposition}
\newtheorem{lemma}[theorem]{Lemma}

\newtheorem{definition}{Definition}
\theoremstyle{definition}
\newtheorem{remark}[theorem]{Remark}

\DeclareMathOperator{\length}{length}

\title[Entropy stability for Bowen-Series-like maps]{Entropy stability and Milnor-Thurston invariants for Bowen-Series-like maps}
\author{Llu\'{\i}s Alsed\`a$^{1,2}$, David Juher$^{3}$, J\'er\^ome Los$^4$ \and Francesc
Ma\~{n}osas$^1$}
\address{$^1$Departament de Matem\`atiques, Edifici C, 08193 Campus
de Bellaterra, Barcelona, Spain}
\email{lluisalsedaisoler@mat.uab.cat, lluisalsedaisoler@crm.cat}
\email{manyosas@mat.uab.cat}

\address{$^2$ Centre de Recerca Matemàtica, Edifici C, 08193 Campus
de Bellaterra, Barcelona, Spain}

\address{$^3$Departament d'Inform\`atica i Matem\`atica Aplicada,
Universitat de Girona, c/ de la Universitat de Girona, 6, 17003 Girona, Spain}
\email{david.juher@udg.edu}

\address{$^4$Aix-Marseille Universit\'e, CNRS,
Institut Mathematiques de Marseille UMR 7373,
39 Rue F. Joliot Curie, 13013  Marseille,
France}
\email{jerome.los@univ-amu.fr}

\thanks{This work has been possible due to the financial support of our respective institutions,
the MICINN grant number PID2020-118281GB-C31 and the Severo Ochoa and María de Maeztu Program
for Centers and Units of Excellence in R\&D (CEX2020-001084-M) funded by MCIN/AEI/10.13039/501100011033.
We thank CERCA Programme/Generalitat de Catalunya for institutional support.}


\begin{document}
\begin{abstract}
We define a family of discontinuous maps on the circle, called Bowen-Series-like maps, for geometric presentations of surface groups. The family has $2N$ parameters, where $2N$ is the number of generators of the presentation. We prove that all maps in the family have the same topological entropy, which coincides with the volume entropy of the group presentation. This approach allows a simple algorithmic computation of the volume entropy from the presentation only, using the Milnor-Thurston theory for one dimensional maps.
\end{abstract}

\maketitle

\section{Introduction}\label{intro}
Let $\Sigma$ be a closed compact %hyperbolic
surface of rank larger than 2 and let $P=\langle X |  R\rangle$ be a presentation
of its fundamental group $G:=\pi_1(\Sigma)$. Since the rank is larger than 2, the surface $\Sigma$ is hyperbolic in the geometrical sense, $G$ is a hyperbolic group in the sense of Gromov \cite{Gr,GdlH} and its boundary $\partial G$ is homeomorphic to the circle $S^1$. We consider the Cayley graph of the group presentation $\textrm{Cay}^1(G,P)$ and the Cayley 2-complex $\textrm{Cay}^2(G,P)$.
The presentation $P$ is called {\em geometric} if $\textrm{Cay}^2(G,P)$ is homeomorphic to a plane. In particular, $\textrm{Cay}^1(G,P)$ is a planar graph. This property is satisfied by the classical presentations of any surface group (see for instance \cite{Sti}).

For a hyperbolic group $G$ with a presentation $P=\langle X | R\rangle$, let $\mathcal{X}$ be the \emph{generating set}, defined as the set of generators of the presentation and their inverses (by abuse of language, sometimes we will use the term \emph{generator} to refer to any element of $\mathcal{X}$). For any $g\in G$, the \emph{length of $g$}, denoted by $\length(g)$, is the number of symbols of a minimal word in the alphabet $\mathcal{X}$ representing $g$. It coincides with the number of edges in any geodesic segment (shortest path) in the Cayley graph $\textrm{Cay}^1(G,P)$ connecting the identity element of $G$ to the vertex $g$. The \emph{growth function}
\begin{equation}\label{sigma-n}
\sigma_m := \operatorname{Card}\set{g \in G}{\operatorname{length}(g) = m},
\end{equation}
which is also the number of vertices at distance $m$ from the identity in the Cayley graph, plays a central role in geometric group theory \cite{CW,FP,DlH}. Its exponential growth rate is called the \emph{volume entropy}, defined as
\begin{equation}\label{volentropy}
h_{\textrm{vol}}(G,P)=\lim_{m \rightarrow \infty} \frac{1}{m}\log(\sigma_m ).
\end{equation}

On the other hand, a geometric presentation $P$ of a surface group $G$ allows to define a dynamical system given by a {\em piecewise homeomorphism} of the circle $\Phi_P:S^1\rightarrow S^1$, in the sense that $S^1$ has a finite partition by intervals and $\Phi_P$, restricted to each interval of the partition, is a homeomorphism onto its image \cite{L}. The construction of the map $\Phi_P$ from a geometric presentation $P$ is based on an idea initially due to Bowen \cite{Bo2} and Bowen-Series \cite{BS}. The dynamical complexity of this map can be measured by its \emph{topological entropy}, a notion first introduced in \cite{AKM,Bo1} for continuous maps of compact metric spaces that can be extended to piecewise monotone and discontinuous maps of the circle \cite{MZ}. The topological entropy of a map $f$ will be denoted by $h_{top}(f)$.

The maps defined in \cite{Bo2,BS} and in \cite{L} satisfy several important properties:
\begin{itemize}
\item The map is Markov and expanding.
\item The map and the group $G$ are orbit equivalent. For this statement the group $G$ is viewed as acting on its Gromov boundary $\partial G = S^1$ by homeomorphisms and is considered as a discrete subgroup of
$\textrm{Homeo}(S^1)$.
\item For the particular map $\Phi_P$ defined in \cite{L}, the topological entropy $h_{top}(\Phi_P)$ of the map and the volume entropy $h_{vol}(G,P)$ of the group presentation are equal.
\end{itemize}

In this paper we define, for a given geometric presentation $P$, a family of \emph{Bowen-Series-like maps} $\Phi_\Theta$ indexed by a set of $2N$ parameters $\Theta$, where $N$ is the number of generators in $P$. The map $\Phi_P$ above is a particular member of this much wider family. We prove then that the topological entropy of any map in the family is equal to the volume entropy (\ref{volentropy}) of the presentation. As we will see, this result has a strong computational implication: the volume entropy of any given geometric presentation $P$ becomes easily computable by a purely combinatorial procedure. In order to define the family of maps $\Phi_\Theta$, let us recall some known facts.

The boundary of a hyperbolic geodesic metric space is a topological, metric space \cite{Gr}. Any point in the boundary
is an equivalence class of geodesic rays that remain at a uniform bounded distance from each others.
For finitely generated groups the metric space is a Cayley graph and the hyperbolic property, as well as the topological boundary, does not depend on the presentation. For our group $G$ with a given presentation, a geodesic ray starting at $\Id_G$ is described as an infinite word $W$ in the generating set $\mathcal{X}$ such that any finite subword of $W$ is geodesic. We denote by $[\zeta]$ a particular expression of a geodesic ray representing the point $\zeta\in\partial G$, and by $[\zeta]_m$ the initial word of length $m\geq 1$ of the infinite word $[\zeta]$.
For $x\in\mathcal{X}$, the {\em cylinder set} (of length one) ${\mathcal{C}}_x $ is the subset of $\partial G$ defined by
\begin{equation}\label{cylinderxx}
{\mathcal{C}}_x  = \lbrace \zeta \in \partial G : \mbox{ there exists } [\zeta]  \mbox{ with } [\zeta]_1 = x \rbrace.
\end{equation}

This notion of cylinder set is naturally extended to all lengths. If $g\in G$ is an element of length $m\geq 1$ and $\{g\}$ is the set of all geodesic expressions of $g$, we denote:
\begin{equation}\label{cylinderx}
{\mathcal{C}}_g  = \lbrace \zeta \in \partial G : \mbox{ there exists } [\zeta]  \mbox{ with }[\zeta]_m \in \{ g \}  \rbrace.
\end{equation}

In combinatorial dynamics a cylinder set, say of one letter, is the set of infinite words starting with that letter.
Here the infinite word is replaced by the notion of geodesic ray and the set of infinite words is replaced by the boundary of the hyperbolic space.
The main difference, if the space is not a tree, is that many geodesic rays might define the same point on the boundary. Therefore cylinder sets intersect in general.
This notion of cylinder sets for hyperbolic spaces is sometimes called ``shadow'' in the geometry literature \cite{Cal}.

In our particular situation of a Cayley graph for a geometric presentation of a surface group, the cylinder sets satisfy particular properties \cite{L}:
\begin{itemize}
\item[(I)] ${\mathcal{C}}_x$ is a connected interval of $S^1=\partial G$ for each generator $x\in\mathcal{X}$.
\item[(II)] ${\mathcal{C}}_x \cap {\mathcal{C}}_y \neq \emptyset$ if and only if $x$ and $y$ are two adjacent generators, with respect to the cyclic ordering induced by the planarity of the Cayley graph
$\textrm{Cay}^1(G,P)$. In addition, if this intersection is non empty, then it is a connected interval.
\end{itemize}

It is well known that a hyperbolic group acts on its boundary by homeomorphisms. Here we obtain that the group $G$ can be seen as a discrete subgroup of $\textrm{Homeo}(S^1)$.

\begin{definition}\label{standing}
Let $\Sigma$ be a closed, compact surface of negative Euler characteristic and let $P$ be a geometric presentation of $G:=\pi_1(\Sigma)$.
We denote the elements of the generating set $\mathcal{X}$ as $x_1,x_2,\ldots,x_{2N}$, where the indices are defined modulo $2N$ (with $1,2,\ldots,2N$ as representatives of the classes modulo $2N$) so that $x_j$ is adjacent to $x_{j\pm 1}$. From the point (II) above there are $2N$ disjoint intervals
\[  J_j := {\mathcal{C}}_{x_{j-1}} \cap {\mathcal{C}}_{x_{j}} \subset S^1. \]
For each $\Theta:=(\theta_1,\theta_2,\ldots,\theta_{2N}) \in  J_1 \times J_2 \times\ldots\times J_{2N}$ we consider the intervals $I_j:=[\theta_{j},\theta_{j+1})\subset S^1$ and the map
\begin{equation}\label{family}
\map{\Phi_{\Theta}}{S^1} \mbox{ so that } \Phi_{\Theta}(z)=x_j^{-1}(z) \mbox{ if } z\in I_j.
\end{equation}
Such a map is called a Bowen-Series-like map. Each point $\theta_i$ will be called a cutting point, and $\Theta$ will be called a cutting parameter.
\end{definition}

The notation $x_j^{-1}(z)$ in Definition~\ref{standing} stands for the action, by homeomorphism, of the group element $x_j^{-1}$ on the boundary $\partial G=S^1$. The map $\Phi_{\Theta}$ is, therefore, a {\em piecewise homeomorphism}. More precisely, the intervals $I_j$ define a partition of the space $S^1$ and $\Phi_{\Theta}\evalat{I_j}$ is a homeomorphism onto its image.

With all these ingredients we are ready to state the main result of this paper.

\begin{MainTheorem}\label{mainTh}
Let $\Sigma$ be a closed, compact surface of negative Euler characteristic and let $P$ be any geometric presentation of $G:=\pi_1(\Sigma)$. Then, for any cutting parameter $\Theta\in J_1\times J_2\times\ldots\times J_{2N}$,
the Bowen-Series-like map $\Phi_{\Theta}$ satisfies:
\begin{enumerate}
\item[(a)] $h_{top}(\Phi_{\Theta})=h_{vol}(G,P)=\log(\lambda)$, where $1/\lambda$ is the smallest root in $(0,1)$ of an integer polynomial $Q_P(t)$ that can be explicitly computed from $P$.
\item[(b)] $\Phi_{\Theta}$ is topologically conjugate to a piecewise affine map $\widetilde{\Phi_{\Theta}}$ whose slope is constant for all intervals $I_j$ and is equal to $\pm \lambda$.
\end{enumerate}
\end{MainTheorem}

A priori, Theorem~\ref{mainTh} is a surprising result, since the dynamics of two different maps in the family are quite different, in particular they are not pairwise topologically conjugate or even semi-conjugate. For some choices of the parameters $\Theta$, the map $\Phi_{\Theta}$ satisfies the Markov property, while for some other choices the map is not Markov.

In addition, for many geometric presentations (in particular for the classical ones), and for an open set of parameters $\Theta$, the corresponding maps
$\Phi_\Theta$ satisfy the assumptions of \cite{LV}, where it is proved that such maps are orbit equivalent to the surface group action on $S^1$. This implies that for many presentations and many parameters $\Theta$, two different maps are orbit equivalent to each other since they are both orbit equivalent to the same group action. The orbit equivalence is a much weaker relation than the conjugacy and, a priori, does not preserve the topological entropy. These observations imply that the family $\Phi_{\Theta}$ shows some surprising stability properties. The question of the orbit equivalence to the group action for all parameters and all geometric presentations is not considered here. Some works are in progress on that problem.

\begin{remark}
The map $\Phi_{\Theta}$ is defined by the action of the generators of the group $G$ on some intervals of the circle, which, as the Gromov boundary of $G$, is just a topological space. After the conjugacy given by Theorem \ref{mainTh}(b), $\widetilde{\Phi_{\Theta}}$ is piecewise affine with a constant slope, the algebraic number $\lambda$. For this new map, the circle, i.e. the boundary of the group, admits a well defined metric that reflects the growth property of the group presentation. This is an intriguing property.
\end{remark}

An interesting part of Theorem~\ref{mainTh} is the computational property stated in (a). Since the polynomial invariant of the presentation $Q_P(t)$ does not depend on the cutting parameter $\Theta$, it is possible to choose a particular $\Theta$ for which the Milnor-Thurston invariants \cite{MT} can be easily computed using elementary algebraic operations. The computation is at the end quite ``simple'' and can be arranged in the form of a purely combinatorial algorithm that takes as input the presentation of the group and prints out the polynomial $Q_P(t)$ and the corresponding volume entropy. The explicit algorithm and some examples are provided in Section~\ref{computer}.

The paper is organized as follows. In Section~\ref{review} we gather the relevant properties of the family $\Phi_{\Theta}$ by reviewing some known facts from previous works. In Section~\ref{treelike} we describe a new ``tree-like'' structure associated to each intersection interval $J_i$. This structure defines some particular cutting parameters $\theta_j\in J_j$ where the map exhibits possibly different dynamical behaviors. Section~\ref{staircases} is concerned with the geometric description of the set of all geodesics connecting two vertices of the Cayley graph. In Section~\ref{central} we define, for each interval $J_i$, an open subinterval $C(J_i)\subset J_i$ that we call \emph{central}. We show that when the cutting parameter $\Theta$ is central, i.e. $\theta_j\in C(J_j)$ for all $j$, the dynamics of $\Phi_{\Theta}$ at each cutting point is specially simple: for some $k_j$, the $k_j$-th iterates of $\Phi_\Theta$ from the left and the right coincide at $\theta_j$. In Section~\ref{middle} we define a particular central parameter which we call the \emph{middle parameter}. In Section~\ref{inequalities} we prove, independently of the parameter $\Theta$, two inequalities comparing the rates of increasing of the two entropy functions, one for the group presentation and one for the dynamics of $\Phi_{\Theta}$. This step, together with the use of the middle parameter and the techniques from the Milnor-Thurston kneading theory, that are summarized in Section~\ref{milnor}, leads to the proof of Theorem~\ref{mainTh}. Finally, in Section~\ref{computer} we fully describe the algorithm that computes the polynomial invariant $Q_P(t)$ and give some explicit examples for several presentations.

\section{Review of the Bowen-Series-like maps $\Phi_\Theta$}\label{review}
In this Section we gather some properties, from \cite{L}, of geometric presentations for hyperbolic surface groups and of the maps of type $\Phi_{\Theta}$. Recall that a presentation
$P=\langle X|R\rangle$ of a surface group is called {\em geometric} if the Cayley 2-complex $\textrm{Cay}^2(G,P)$ is planar. The next result (Lemma~2.1 of \cite{L}) states some elementary consequences of this planarity property.

\begin{lemma}\label{co}
Let $G$ be a co-compact hyperbolic surface group and let $P=\langle X|R\rangle=\langle g_1,\ldots,g_N|R_1,\ldots,R_k\rangle$ be a geometric presentation of $G$. The following conditions hold.
\begin{enumerate}
\item[(a)] The set $\{g_1^{\pm1},g_2^{\pm1},\ldots,g_N^{\pm1}\}$ admits a cyclic ordering that is preserved by the group action.
\item[(b)] There exists a planar fundamental domain $\Delta_P$.
\item[(c)] Each generator appears exactly twice (with $+$ or $-$ exponent) in the set of relations $R$.
\item[(d)] Let $a,b$ be a pair of adjacent generators according to the cyclic ordering given by (a). Then, there is exactly one relation $R_i$ such that a cyclic shift of $R_i$ contains either $b^{-1}a$ or $a^{-1}b$ as a sub-word.
\end{enumerate}
\end{lemma}

\begin{figure}
\centering
\includegraphics[scale=0.7]{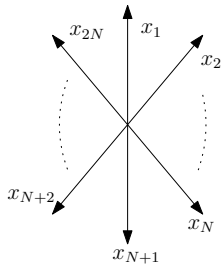}
\caption{The clockwise labelling of the generating set associated to the cyclic ordering given by
Lemma~\ref{co}(a).}
\label{fig:notationCicOrd}
\end{figure}

Here we recall the standing notation introduced in Definition~\ref{standing}. The elements of the generating set $\mathcal{X}:=\{g_1^{\pm1},g_2^{\pm1},\ldots,g_N^{\pm1}\}$ are denoted as $x_1,x_2,\dots,x_{2N}$ in such a way that $x_{i\pm1}$ are the elements adjacent to $x_i$ with respect to the cyclic ordering from Lemma~\ref{co}(a). We also adopt the convention that $x_i$ is on the \emph{left} of $x_{i+1}$ (see Figure~\ref{fig:notationCicOrd}). This convention defines a clockwise orientation of the plane $\mathrm{Cay}^2(G,P)$.

It is easy to see that a presentation containing a relation of the form $x_i^2$ does not satisfy Lemma~\ref{co}(a), and thus it cannot be geometric. On the
other hand, a relation of the form $x_ix_j$ is simply a trivial identification of the generators $x_i$ and $x_j^{-1}$. So, all presentations considered
in this paper will not contain relations of length two by hypothesis. The following purely technical observation will be necessary to exclude unwanted
particular situations in some proofs.

\begin{remark}\label{restriction}
Assume that our presentation $P$ has $N=3$ generators. By Lemma~\ref{co}(c), the sum of the lengths of all relations in $P$ is 6. It is easy to see that
a presentation with 3 generators and 2 relations of length 3 corresponds either to a torus or to a Klein bottle, the rank 2 cases, that are not hyperbolic.
Since we are excluding by principle relations of length 2, it follows that the only possible geometric presentation with $N=3$ generators has only
one relation of length 6.
\end{remark}

A {\em bigon} in $\textrm{Cay}^1(G,P)$ between two vertices $v,w$ is a pair $\{\gamma,\gamma'\}$ of disjoint geodesic segments connecting $v$ to $w$.
% More precisely, $\gamma$ and $\gamma'$ are words of minimal length in the alphabet $\mathcal{X}$ such that $v\gamma=v\gamma'=w$ and $\gamma$ and $\gamma'$ have no letters in common.
We denote by $B(x,y)$ the set of bigons, starting at $\Id$, so that one of the two geodesics starts with $x$ and the other one starts with $y$. Note that both words have the same length, which will be called the \emph{length} of the bigon. A \emph{bigon ray} is a pair of disjoint geodesic rays connecting a vertex in $\textrm{Cay}^1(G,P)$ to a point on the boundary $\partial G$.
The following result is a combination of Lemma~2.8 and Corollary~2.11 of \cite{L}.

\begin{proposition}\label{Bigons}
Let $P=\langle X|R\rangle$ be a geometric presentation of a hyperbolic surface group $G$. Then,
\begin{enumerate}
\item [(a)] $B(x,y)\neq\emptyset$ if and only if $x$ and $y$ are adjacent in the cyclic ordering of $\mathcal{X}$.
\item[(b)] If $x$ and $y$ are adjacent generators in $\mathcal{X}$ then there exists a unique bigon $\beta(x,y)\in B(x,y)$ of minimal length.
\end{enumerate}
\end{proposition}

A point $z\in\partial G$ is the limit of possibly many geodesic rays starting at the identity. We denote by $[z]$ the expression of one geodesic ray in $\textrm{Cay}^1(G,P)$ converging to $z$, this is an infinite word in the alphabet $\mathcal{X}$.

Recall from (\ref{cylinderxx}) that the {\em cylinder set} ${\mathcal{C}}_x$, for a generator $x\in\mathcal{X}$, is the subset of points $\zeta\in\partial G$ such that there exists $[\zeta]$, a ray at $\Id$ converging to $\zeta$, such that $[\zeta]$ starts with $x$.

\begin{figure}
\centering
\includegraphics[scale=0.7]{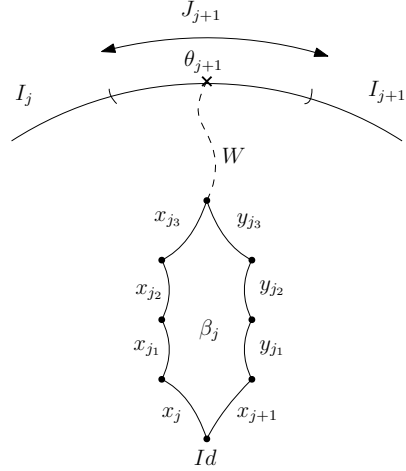}
\caption{The bigon $\beta(x_j,x_{j+1})$ and the cutting point $\theta_{j+1}$.}
\label{fig:1}
\end{figure}

The following result is a direct consequence of Section~3 of \cite{L}. It is illustrated by Figure~\ref{fig:1}.

\begin{proposition}\label{cylinders}
The cylinder sets satisfy:
\begin{enumerate}
\item[(a)] For any $x\in\mathcal{X}$, ${\mathcal{C}}_x$ is connected and ${\mathcal{C}}_x \cap {\mathcal{C}}_y\neq\emptyset$  if and only if $x$ and $y$ are adjacent generators. In this case it is an interval.
\item[(b)] For any $\theta\in {\mathcal{C}}_x \cap {\mathcal{C}}_y$, there is an infinite word $W$ in the alphabet $\mathcal{X}$ so that $\theta \in \partial G$ has two geodesic ray expressions
$[\theta]_{-}=L_{x}W$ and $ [\theta]_{+}=L_{y}W$, where $ \{ L_{x}, L_{y} \}$  are the two geodesic segments of the minimal bigon $\beta(x,y)$.
\end{enumerate}
\end{proposition}

In Proposition~\ref{cylinders}, the infinite word $W$ might not be unique and the possible non uniqueness will play a key role later.

The map $\Phi_{\Theta}$ defined in (\ref{family}) is a piecewise homeomorphism of $S^1=\partial G$. Each interval $I_j$ of the partition of $S^1=\bigcup_{j=1}^{2N} I_j$ is defined by $I_j=[\theta_{j},\theta_{j+1})$, where
$\theta_j \in {\mathcal{C}}_{x_{j-1}} \cap {\mathcal{C}}_{x_{j}} $ is called a {\em cutting point}. See Figure~\ref{km1} for a particular example. The map is defined as
\[ \Phi_{\Theta} (z) = x_j^{-1} (z) \mbox{ if } z \in I_j. \]

The definition implies, for each $j\in\{1,2,\ldots,2N\}$, that $I_j \subset {\mathcal{C}}_{x_j}$ and the map
${\Phi_{\Theta}}\evalat{I_j}$ is a homeomorphism onto its image. At the cutting points the map is not continuous.
From the definition of the cylinder sets ${\mathcal{C}}_{x_j}$, if $z\in I_j$ then there exists a geodesic ray, starting at $\Id$ in the Cayley graph, denoted by $[z]$ and converging to $z\in\partial G$, such that
$[z]=x_jA$, where $A$ is an infinite word in the alphabet $\mathcal{X}$.
The map applied to the point $z$ thus gives:
\begin{equation}\label{shiftMap}
\textrm{if } z \in I_j  \textrm{ then } [z]=x_jA \textrm{  and  } [ \Phi_{\Theta}(z)] = A.
\end{equation}
In other words $\Phi_{\Theta}$ is a standard ``shift map'' for the particular coding of the boundary $\partial G$ obtained from the intervals $I_j$.
Different choices of the cutting parameter $\Theta$ in $J_1\times J_2\times\ldots\times J_{2N}$ define different maps and thus different codings of $\partial G$.

There is another important consequence of the results in \cite{L} that will be used to prove statement (b) of Theorem~\ref{mainTh}. As it has been mentioned in
Section~\ref{intro}, the map $\Phi_P$ defined in \cite{L} is nothing but a Bowen-Series-like map $\Phi_\Theta$ for a particular choice of the parameter $\Theta$.
Among other properties, Proposition~5.3 of \cite{L} states that, when no generator appears twice in a relation of length 3, $\Phi_P$ is topologically transitive
in the sense that, for any two open intervals $U,V\subset S^1$, there exists $n\ge0$ such that $\Phi_P^n(U)\cap V\ne\emptyset$. It is easy to see that the
arguments for proving the transitivity can be directly generalized to any $\Phi_\Theta$ independently of the parameter $\Theta$, and that the condition
about the generators is, from the point of view of the transitivity, also irrelevant. So, we have the following result.

\begin{proposition}\label{transitive}
For any parameter $\Theta$, the Bowen-Series-like map $\Phi_\Theta$ is topologically transitive.
\end{proposition}

\section{A tree-like structure for each parameter interval}\label{treelike}
In this Section we study some particular sequences of points in the intervals
$J_j = {\mathcal{C}}_{x_{j-1}} \cap {\mathcal{C}}_{x_{j}}$. These points are limits of a tree-like set in the Cayley graph. They are interpreted as the parameters where the dynamics of the map could potentially change, so in some sense they are like ``bifurcation'' parameters. The starting point is a characterization of the extreme points of $J_j$ obtained in \cite{BD}. In order to simplify the notations and avoid too many indices we consider a pair of adjacent generators $(x,y)$ and the corresponding minimal bigon $\beta(x,y)$ with the convention that the generator $x$ is on the left of $y$ in the cyclic ordering induced by our clockwise orientation of $S^1$. We consider also the minimal bigons $\beta_v(x,y)$ starting at a vertex $v$, possibly different from the identity, with length denoted by $k(x,y)$. To fit with these notations we denote $J(x,y)$ the interval ${\mathcal{C}}_{x} \cap {\mathcal{C}}_{y}$.

Assume that the minimal bigon $\beta (x, y ) = \{ L_{x} , L_{y} \}$ given by Proposition~\ref{Bigons}(b) is:
\begin{equation}\label{bigons}
 L_{x} = x\cdot x_{j_2}\cdots x_{j_k}  \textrm{ and }
 L_{y} = y\cdot y_{j_2}\cdots y_{j_k} \textrm{, with } k = k(x,y),
\end{equation}
where the various $x_{j_m},y_{j_n}$ are generators in $\mathcal{X}$.

Consider the last edge along the geodesic segment $L_x$, with label $x_{j_k}$. It starts at the vertex $v_{k-1}^{0}$, which is also the terminal vertex of the geodesic path labelled $x\cdot x_{j_2}\cdots x_{j_{k-1}}$ starting at the identity. Symmetrically consider the last edge, labeled $y_{j_k}$, starting at the vertex $v_{k-1}^{1}$ along the segment $L_y$.

Consider the minimal bigon $\beta_{v_{k-1}^{0}}(x_{j_1}^{0},x_{j_k})$, based at $v_{k-1}^{0}$, where we denote $x_{j_1}^{0}$ the generator that is adjacent to $x_{j_k}$ on the left at the vertex $v_{k-1}^{0}$. Symmetrically we consider $\beta_{v_{k-1}^{1}}(y_{j_k},y_{j_1}^{1})$, the minimal bigon based at $v_{k-1}^{1}$, and the generator $y_{j_1}^{1}$ adjacent to $y_{j_k}$ on the right.

We define now two binary operations on bigons that are called {\em extensions}: \label{exten}
\begin{eqnarray}\label{Extension}
 {\mathcal{E}} (0; \beta(x,y) ) := \beta (x,y)  \otimes_{x_{j_k}}
\beta_{ v_{k-1}^{0} }( x_{j_1}^{0}, x_{j_k}), \nonumber
\\
 {\mathcal{E}} (1; \beta(x,y) ) := \beta (x,y) \otimes_{y_{j_k}}
\beta_{ v_{k-1}^{1} }( y_{j_k}, y_{j_1}^{1}),
\end{eqnarray}
where the symbol $\otimes_{a}$ means the concatenation of the two bigons along their common edge $a$. More precisely, if $\beta_{ v_{k-1}^{0} }( x_{j_1}^{0}, x_{j_k})=\{x_{j_1}^{0} x_{j_2}^{0}\cdots x_{j_{k^0}}^{0};x_{j_k} y_{j_2}^{0}\cdots y_{j_{k^0}}^{0}\}$, the 0-extension is:
\begin{equation}\label{extension0}
{\mathcal{E}} (0; \beta(x,y) ) = \{ x\cdot x_{j_2}\cdots x_{j_{k-1}}\ast  x_{j_1}^{0} x_{j_2}^{0}\cdots x_{j_{k^0}}^{0};y\cdot y_{j_2}\cdots y_{j_k} \ast  y_{j_2}^{0}\cdots y_{j_{k^0}}^{0} \}
\end{equation}
where the symbol $\ast$ represents the location where the paths are concatenated. Symmetrically, if $\beta_{ v_{k-1}^{1} }( y_{j_k}, y_{j_1}^{1})=\{ y_{j_k} x_{j_2}^{1}\cdots x_{j_{k^1}}^{1} ;  y_{j_1}^{1} y_{j_2}^{1}\cdots y_{j_{k^1}}^{1}\}$, the 1-extension is
\begin{equation}\label{extension1}
 {\mathcal{E}} (1; \beta(x,y) ) = \{ x\cdot x_{j_2}\cdots x_{j_{k}}\ast  x_{j_2}^{1}\cdots x_{j_{k^1}}^{1}; y\cdot y_{j_2}\cdots y_{j_{k-1}} \ast  y_{j_1}^{1}\cdots y_{j_{k^1}}^{1} \}.
\end{equation}

\begin{figure}%[htbp]
\centering
\includegraphics[scale=0.65]{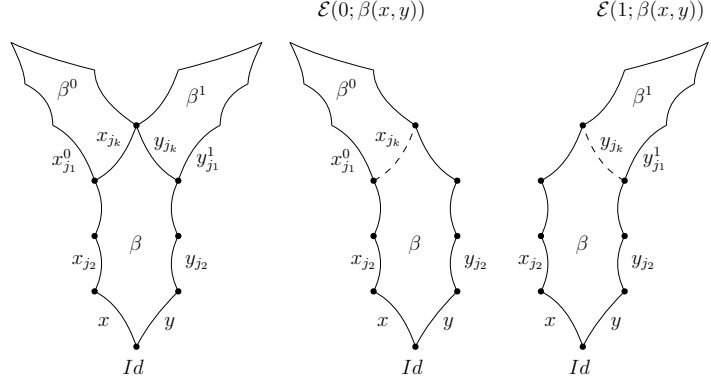}
\caption{The left (0) and right (1) extensions of a minimal bigon.}
\label{fig:2}
\end{figure}

These operations are represented in Figure~\ref{fig:2}. Observe that the minimal bigons occurring above have, a priori, different lengths, denoted as $k=k(x,y)$, $k^0=k(x_{j_1}^{0},x_{j_k})$ and $k^1=k(y_{j_k},y_{j_1}^{1})$.

\begin{proposition}\label{extension}
The extensions $ {\mathcal{E}} (0; \beta(x,y) )$ and $ {\mathcal{E}} (1; \beta(x,y) )$ are well defined bigons in $B(x,y)$.
\end{proposition}
\begin{proof}
We first observe that the four concatenations of paths in (\ref{extension0}) and (\ref{extension1})
are well defined, i.e. the terminal vertex of the initial path coincides with the initial vertex of the next one. The four concatenations of (\ref{extension0}) and (\ref{extension1}) start at the identity. For the 0-extension, the terminal vertex of the two paths in (\ref{extension0}) is the terminal vertex of the bigon
$\beta_{ v_{k-1}^{0} }( x_{j_1}^{0}, x_{j_k})$ and similarly the terminal vertex of the 1-extension in (\ref{extension1}) is the terminal vertex of $\beta_{v_{k-1}^{1}}(y_{j_k}, y_{j_1}^{1})$.

It remains to check that each concatenation is a geodesic. This is a direct consequence of \cite{L} (see Lemma~2.9), more details can be found in \cite{BD}. The two extensions are thus two bigons in the set $B(x,y)$.
\end{proof}

To simplify the notations we now write the two extensions as ${\mathcal{E}} (0/1; \beta )$. The extension operations start with a bigon $\beta\in B(x,y)$ and give two well defined bigons in $B(x,y)$.
Therefore we can iterate the extension construction.

We simplify again the notations and we denote
\[ {\mathcal{E}} (\epsilon_2, \epsilon_1 ; \beta ) :=
{\mathcal{E}} ( \epsilon_2;  {\mathcal{E}} (\epsilon_1 ; \beta ) ) \mbox{ for }
(\epsilon_1, \epsilon_2) \in \{0, 1\}^2, \]
and by a finite iteration we obtain
\begin{equation}\label{extensionN}
{\mathcal{E}} (\epsilon_n,\ldots,  \epsilon_1 ; \beta ) :=
{\mathcal{E}} ( \epsilon_n;  {\mathcal{E}} (\epsilon_{n-1},\ldots,  \epsilon_1 ; \beta ) ) \textrm{ for }
(\epsilon_1, \ldots, \epsilon_n) \in \{0, 1\}^n.
\end{equation}

\begin{lemma}\label{limitExtension}
Let $(x,y)$ be a pair of adjacent generators. Then,
\begin{enumerate}
\item[(a)] For all $n \in \mathbb{N}$ and all $(\epsilon_n,\ldots,\epsilon_1)\in\{0, 1\}^n$, the collection of level $n$ extensions
${\mathcal{E}} (\epsilon_n,\ldots,\epsilon_1;\beta(x,y))$ are well defined bigons in $B(x,y)$.
\item[(b)] The limit sequences, when $n\rightarrow \infty$, define a collection of bigon rays in $B(x,y)$ and, for each
${\bf E }\in  \{0, 1\}^{\mathbb{N}}$, the limit is a well defined point
\[ \Omega ({\bf E } ; \beta(x,y) ) \in {\mathcal{C}}_x \cap {\mathcal{C}}_y \subset \partial G. \]
\end{enumerate}
\end{lemma}
\begin{proof}
By Proposition \ref{Bigons} there is a unique minimal bigon for each pair of adjacent generators. There are thus finitely many minimal bigons, each one with a length $k(j) \geq 2$.
Let us prove (a). For each $n \in {\mathbb{N}}$ and each $(\epsilon_n, \ldots, \epsilon_1) \in \{0, 1\}^n$ the finite extension ${\mathcal{E}}(\epsilon_n,\ldots, \epsilon_1 ; \beta(x,y))$ is well defined by repeatedly
applying Proposition~\ref{extension}. The length of each such bigon is an explicit additive function of the lengths of the various bigons $\beta_v$ occurring in the finite extension
${\mathcal{E}}(\epsilon_n,\ldots,\epsilon_1;\beta(x,y))$. For instance in (\ref{Extension}) the length is $k(x,y)+k(x_{j_1}^{0},x_{j_k})-1$. The length of the bigons
${\mathcal{E}} (\epsilon_n,\ldots,\epsilon_1;\beta(x,y))$ is strictly increasing with $n$, since $k(j)\geq 2$ for each $j$.

To prove (b), note that ${\mathcal{E}}(\epsilon_n,\ldots,\epsilon_1;\beta)$ define sequences of bigons in $B(x,y)$ whose lengths go to infinity when $n\rightarrow \infty$. For each ${\bf{ E }} \in  \{0, 1\}^{\mathbb{N}}$, the set of geodesic rays in ${\mathcal{E}} ({\bf{ E }} ; \beta(x,y) )$ remains at a finite distance from each other. More precisely let $K:=\max\big\{k(m) ; m\in \{1,2,\ldots,2N\}\big\}$, where $k(m)$ is the length of a minimal bigon. Then two points, at distance $M$ from the identity along the geodesic rays of ${\mathcal{E}} ({\bf{ E }} ; \beta(x,y))$, remain at distance less than $K$ from each other, for all $M$. By definition of the Gromov boundary \cite{Gr,GdlH} the geodesic rays in ${\mathcal{E}} ({\bf{ E }} ; \beta(x,y) )$, for each ${\bf{ E}} \in  \{0, 1\}^{\mathbb{N}}$, converge to the same point in $\partial G$. In addition, each ${\mathcal{E}} (\epsilon_n,\ldots,\epsilon_1;\beta(x,y))$ is a bigon in $B(x,y)$ and, by definition of the cylinder sets, the limit point  $\Omega({\bf{E}};\beta(x,y))$ belongs to ${\mathcal{C}}_x\cap{\mathcal{C}}_y\subset\partial G$.
\end{proof}

As it is clear from the definition, the infinite extension construction is a ``tree-like'' process and the underlying tree is rooted and binary.

\begin{lemma}\label{orderPreserving}
For any adjacent pair of generators $(x,y)$ there is a well defined map
\[ \begin{array}{rcl}
\Omega :  \big\{0, 1\big\}^{\mathbb{N}} & \longrightarrow & {\mathcal{C}}_x \cap {\mathcal{C}}_y \subset \partial G \\
{\bf{ E}} \in  \{0, 1\}^{\mathbb{N}} & \longrightarrow & \Omega ({\bf{ E}} ; \beta(x,y))
\end{array} \]
satisfying:
\begin{enumerate}
\item[(a)] $\Omega$ is injective
\item[(b)] $\Omega$ is order preserving on $\partial G$ for the ordering of $\{0, 1\}^{\mathbb{N}}$ induced by $ 0 < 1$.
\end{enumerate}
\end{lemma}
\begin{proof} The fact that $\Omega$ is a well defined map is a direct consequence of Lemma~\ref{limitExtension}. Let us prove (a). The fact that $\Omega$ is injective comes from Proposition~\ref{Bigons}(a). Indeed, if
${\bf{E}} \neq {\bf{ E ' }}$ in $\big\{0, 1\big\}^{\mathbb{N}}$ then there is a first integer $n$ so that $\epsilon_n \neq \epsilon'_n$. From the extension construction, there is a vertex $v_n$ in the Cayley graph such that
$v_n \in {\mathcal{E}} ( \epsilon_n, \ldots, \epsilon_1 ; \beta(x,y) ) \cap
{\mathcal{E}} ( \epsilon'_n, \ldots, \epsilon'_1 ; \beta(x,y) )$ and, starting from $v_n$, the geodesic rays, say $\gamma$ and $\gamma'$, are different in ${\mathcal{E}}({\bf{E}};\beta(x,y))$ and in
${\mathcal{E}}({\bf{E'}};\beta(x,y))$. Here being different means that the two geodesic rays based at $v_n$ start by two non adjacent generators (see Figure~\ref{fig:2} for $n=1$). These two geodesic rays cannot converge to the same point in $\partial G$ since otherwise they would define a bigon ray, which is impossible by Proposition~\ref{Bigons}(a).

Let us prove (b). If ${\bf{E}}<{\bf{E'}}$ then, by the arguments above, there is $n$ such that $\epsilon_n<\epsilon'_n$, meaning $\epsilon_n=0$ and $\epsilon'_n=1$. The two geodesic rays
$\gamma$ and $\gamma'$ described above and based at the vertex $v_n$ are such that
$\gamma$ starts at $v_n$ by a generator $y^{0\ast}_{j_2}$ and $\gamma'$ starts at $v_n$ by a generator $x^{1\ast}_{j_2}$, with the notations (\ref{extension0}) and (\ref{extension1}) (for $n=1$, see Figure~\ref{fig:2}). By the cyclic ordering at $v_n$ and with our conventions we observe that $y^{0\ast}_{j_2}<x^{1\ast}_{j_2}$. By Proposition~\ref{Bigons}(a), the two geodesic rays $\gamma$ and $\gamma'$ cannot intersect after the vertex $v_n$ since otherwise they would define a bigon starting with non adjacent generators at $v_n$. Therefore, the two rays converge to two points on $\partial G$ and the cyclic ordering at $v_n$ is preserved on the boundary. Thus we obtain
$\Omega({\bf{E}};\beta(x,y))<\Omega({\bf{E'}};\beta(x,y))$.
\end{proof}

The idea of the extension construction appeared in \cite{BD} for the special cases ${\mathcal{E}}(0^{\infty};\beta)$ and ${\mathcal{E}}(1^{\infty};\beta )$, where they are called ``ladders''.

\begin{theorem}\label{BamDia}
With the above notations, $\Omega(0^{\infty};\beta(x,y))$ and $\Omega(1^{\infty};\beta(x,y))$ are the two extreme points of the interval
$J(x,y)={\mathcal{C}}_x \cap {\mathcal{C}}_y \subset \partial G$.
\end{theorem}

This result from \cite{BD} is a characterization of the intersection of cylinder sets and therefore of the cylinder sets. These special points are thus quite important for the study of the family $\Phi_{\Theta}$.

\section{Staircases and galleries}\label{staircases}
In this Section we analyze the geometric structure of the subgraph of the Cayley graph composed by all geodesics joining two
given elements of $G$. The characterization of this subgraph will be one of the tools to prove statement (a) of Theorem~\ref{mainTh}.

For any pair $v,w\in G$ consider a geodesic segment in the Cayley graph $\textrm{Cay}^1(G,P)$ connecting the vertices $v$ and $w$. We keep the same notations $v,w$ for the corresponding elements in the group $G$. It is represented as a word $W$ of minimal length in the generating set $\mathcal{X}$. The \emph{distance} between $v$ and $w$, denoted as $d(v,w)$, is the number of symbols in $W$.

\begin{figure}
\centering
\includegraphics[scale=0.55]{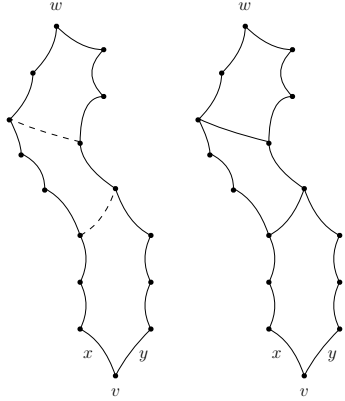}
\caption[fck]{Left: an extension $\mathcal{E}(0,1;\beta(x,y))$ between $v$ and $w$. The dashed edges do not belong to the extension. Right: the corresponding staircase between $v$ and $w$.}\label{stair}
\end{figure}

Let $v,w\in G$. A subgraph $S$ of $\textrm{Cay}^1(G,P)$ will be called a \emph{staircase between $v$ and $w$} if
$S$ is either a minimal bigon between $v$ and $w$ or the union of all minimal bigons concatenated in an \emph{extension} (see Lemma \ref{limitExtension})
between $v$ and $w$. See Figure~\ref{stair} for an example. Note that
if $S_{v,w}$ is a staircase between $v$ and $w$ and $d(v,w)=n$ then for all $1\le j<n$ there are
exactly two vertices in $S_{v,w}$ at distance $j$ from $v$.

\begin{remark}\label{no3}
Let $L$ and $R$ be the left and right geodesic segments connecting $v$ to $w$ and defining a staircase between $v$ and $w$.
For every vertex $x\in L\cup R$, the number of edges in the Cayley graph departing from $x$ and contained in the bounded component of
$\mathbb{R}^2\setminus(L\cup R)$ is at most two. Having two of such edges is only possible when the presentation $P$
has relations of length 3.
\end{remark}

Let $S_{v,w}$ be a staircase between $v$ and $w$ associated to an extension $\mathcal{E}$ from $v$ to $w$. By Lemma~\ref{limitExtension}(a), $\mathcal{E}$ is a bigon between $v$ and $w$. The left and right
geodesic segments from $v$ to $w$ defining this bigon will be respectively denoted by $S_{v,w}^L$ and $S_{v,w}^R$. Note that in general $S_{v,w}$ contains other geodesics going from $v$ to $w$.

\begin{figure}
\centering
\includegraphics[scale=0.55]{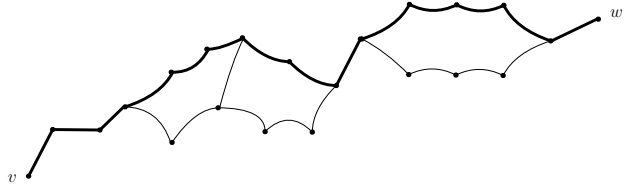}
\caption[fck]{A gallery $F_{v,w}$ between $v$ and $w$. The fat edges correspond to the geodesic path $F^L_{v,w}$.}\label{gallery}
\end{figure}

A {\em gallery} between $v,w\in G$ is a subgraph of the Cayley graph given by a finite set of $k$ staircases $S_{v_1,w_1},\ldots,S_{v_k,w_k}$
and $k+1$ simple paths $\gamma_{v,v_1},\gamma_{w_1,v_2},\ldots,\gamma_{w_k,w}$ (some of them possibly degenerated to a point) such that
\[ \gamma_{v,v_1}S_{v_1,w_1}^L\gamma_{w_1,v_2}\ldots S_{v_k,w_k}^L
\gamma_{w_k,w}\, \textrm{ and } \,\gamma_{v,v_1}S_{v_1,w_1}^R\gamma_{w_1,v_2}\ldots
S_{v_k,w_k}^R \gamma_{w_k,w} \]
are geodesics. See Figure~\ref{gallery} for an example.

If $F_{v,w}$ is a gallery between $v$ and $w$ we denote by $F_{v,w}^L$ and $F_{v,w}^R$
the corresponding geodesic paths defined above and call them the left and right sides of $F_{v,w}$.

Note that if $F_{v,w}$ is a gallery between $v$ and $w$ and
$d(v,w)=n$ then for all $1\le j<n$ there are at most two vertices in
$F_{v,w}$ at distance $j$ from $v$.

\begin{lemma}\label{min3} Let $L$ be a geodesic path between $v,w\in G$. Then, there are at least three different edges
$x_{i_1},x_{i_2},x_{i_3}$ such that $Lx_{i_1},Lx_{i_2},Lx_{i_3}$ are geodesics paths. \end{lemma}
\begin{proof}
This result is a direct consequence of Lemma 2.9 in \cite {L}.
\end{proof}

\begin{lemma}\label{jer1}
Let $L_1,L_2$ be two geodesic paths between $v,w\in G$ such that $L_1\cap L_2=\{v,w\}$. Then, the bounded component of
$\mathbb{R}^2\setminus \{L_1\cup L_2\}$ does not contain vertices of
the Cayley graph.
\end{lemma}
\begin{proof}
Assume, for sake of contradiction, that there exists a
vertex $g$ in the bounded component of $\mathbb{R}^2\setminus \{L_1\cup L_2\}$. Let $W$
be a geodesic path between $v$ and $g$. By Lemma~\ref{min3}, there exist at least three different
edges $x_{i_1},x_{i_2},x_{i_3}$ such that $Wx_{i_1},Wx_{i_2},Wx_{i_3}$ are still geodesic paths.
By applying iteratively Lemma~\ref{min3}, each of these three geodesic paths can be continued up to
a point in $L_1\cup L_2$. From each of these three points we can continue along $L_1$ or $L_2$ up to $w$.
We get then three geodesics from $v$ to $w$ passing through $g$. It follows that we have
three different geodesics from $g$ to $w$ starting respectively with the edges $x_{i_1},x_{i_2},x_{i_3}$.
Two of these edges are not adjacent. In consequence, we have obtained a bigon associated to a pair
of non-adjacent edges, in contradiction with Proposition~\ref{Bigons}.
\end{proof}

\begin{lemma}\label{jer2}
Let $L_1,L_2$ be two geodesic paths between $v,w\in G$ such that $L_1\cap L_2=\{v,w\}.$ Then, $L_1$ and $L_2$ are the left and right sides of a staircase between $v$ and $w$.
\end{lemma}
\begin{proof}
Let $n=d(v,w)$. For all $1\le i<n$, let $v_i$ (respectively $w_i$)
be the point in $L_1$ (resp. $L_2$) satisfying $d(v_i,v)=i$ (resp.
$d(w_i,v)=i$).

\begin{figure}
\centering
\includegraphics[scale=0.7]{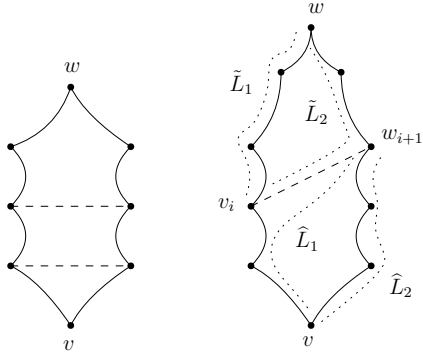}
\caption[fck]{The two cases in the proof of Lemma~\ref{jer2}.}\label{sloping}
\end{figure}

We prove the lemma by induction on $n$. So, assume that the result is true for all lengths less than $n$.
Let $A$ be the bounded component of $\mathbb{R}^2\setminus(L_1\cup L_2)$. By Lemma \ref{jer1}, any
path of the Cayley graph contained in $A$ reduces to a single edge. Since $L_1$
and $L_2$ are geodesics, there are no edges contained in $A$ joining points in $L_1$
(respectively, $L_2$). So, if there exists some edge in $A$, it must connect a point in $L_1$ to a point in
$L_2$. Again, since $L_1$ and $L_2$ are geodesics, any edge that
begins at $v_i$ and ends at $L_2$ must go either to $w_i$, and we call
this edge \emph{horizontal}, or to $w_j$ with $j\in\{i+1,i-1\}$, and we call
this edge \emph{sloping}. If there are no sloping edges, then the pair $\{L_1,L_2\}$ defines
a minimal bigon and there is nothing to prove. See Figure~\ref{sloping}(left). Otherwise,
let $i$ be the smallest index such that there is a sloping edge beginning at $v_i$. Assume
for instance that there is an edge from $v_i$ to $w_{i+1}$ (the other case follows by a
symmetric argument). See Figure~\ref{sloping}(right). Then, the paths $\tilde L_1:=v_{i+1}\cdots v_{n-1}$ and $\tilde
L_2:=w_{i+1}\cdots w_{n-1}$ between $v_i$ and $w$ satisfy the hypothesis of the
lemma and have length less than $n$. By the induction hypothesis, they are the left and right
sides of a staircase. On the other hand, the pair $\{\widehat{L}_1,\widehat{L}_2\}$, where $\widehat{L}_1:=v_1 v_2\cdots
v_i$ and $\widehat{L}_2:=w_1 w_2\cdots w_i$, defines a minimal bigon between $v$ and $w_{i+1}$, since there are
no sloping edges in the interior of the bounded component of $\mathbb{R}^2\setminus(\widehat{L}_1\cup\widehat{L}_2)$.
So, the result follows from the definition of a staircase.
\end{proof}

\begin{lemma}\label{jer3}
Let $L_1,L_2$ be two different geodesic paths between $v,w\in G$. Then, $L_1$ and $L_2$ are the left and right
sides of a gallery between $v$ and $w$.
\end{lemma}
\begin{proof}
We prove the lemma by induction on $n$, the length of the
geodesics. If $n$ is less than the minimum length of a minimal bigon, then
there is nothing to prove. If $n$ is the minimum length of a minimal
bigon, then the pair $\{L_1,L_2\}$ defines a minimal bigon and the
result follows. Assume now that the result is true for all lengths smaller than $n$.
If $L_1\cap L_2=\{v,w\}$, then the result follows from Lemma~\ref {jer3}. If not,
let $z\in L_1\cap L_2$. Then, the result follows by applying the induction hypothesis to
$\tilde L_1$ and $\tilde L_2$, the portions of $L_1$ and $L_2$ joining $v$ and $z$,
and applying also the induction hypothesis to $\widehat{L}_1$ and $\widehat{L}_2$,
the portions of $L_1$ and $L_2$ joining $z$ and $w$.
\end{proof}

\begin{proposition}\label{totes}
The subgraph of all geodesics joining two elements of $G$ is a gallery.
\end{proposition}
\begin{proof}
Let $H$ be the set of all geodesics from $v$ to $w$. Let $L_1$ and $L_2$ be respectively
the leftmost and the rigthmost geodesic paths in $H$ from $v$ to $w$. From
Lemma~\ref{jer3}, they are the left and right sides of some gallery $F_{v,w}$. From the definition
of a gallery, it easily follows that $H=F_{v,w}$.
\end{proof}

The results of this section can be naturally extended to infinite geodesic rays. Recall that a \emph{bigon ray}
has been defined in Section~\ref{review} as a pair of \emph{disjoint} infinite geodesic rays starting (for example)
at the identity element $\Id$ and representing the same point in $S^1$. In Lemma~\ref{limitExtension}(b)
we showed that if $L$ and $R$ are respectively the left and right sides of an infinite extension of bigons, then
the pair $\{L,R\}$ defines a bigon ray. The corresponding point in $S^1$ was denoted by $\Omega({\bf E};\beta(x,y))$
in the statement of Lemma~\ref{limitExtension}(b), where ${\bf E}\in\{0,1\}^{\mathbb{N}}$ is the coding for the
infinite extension and $\beta(x,y)$ is the initial bigon of the extension. The arguments used in the proofs of
Lemmas~\ref{jer1} and \ref{jer2} can be trivially adapted to prove the following converse  result, that states
that \emph{any} bigon ray has precisely this form.

\begin{lemma}\label{jer4}
Let $L,R$ be two disjoint infinite geodesic rays starting at $\Id$. If $L$ and $R$ represent the same point
$\zeta\in S^1$, then $L$ and $R$ are respectively the left and right sides of an infinite extension of
bigons. That is, there exist a pair of adjacent generators $x,y$ and ${\bf E}\in\{0,1\}^{\mathbb{N}}$ such that
$\zeta=\Omega({\bf E};\beta(x,y))$.
\end{lemma}

\section{Dynamical properties in the family $\Phi_{\Theta}$}\label{central}
In this section we study how the dynamics of the maps $\Phi_{\Theta}$ might change with the parameters $\Theta\in J_1\times J_2\times\ldots\times J_{2N}$.  We show in particular the importance of the points in the Cantor set
$\Omega({\bf{E}};\beta(x,y))\in J(x,y)$ as defined in Lemma \ref{limitExtension}.

\begin{figure}%[htbp]
\centering
\includegraphics[scale=0.55]{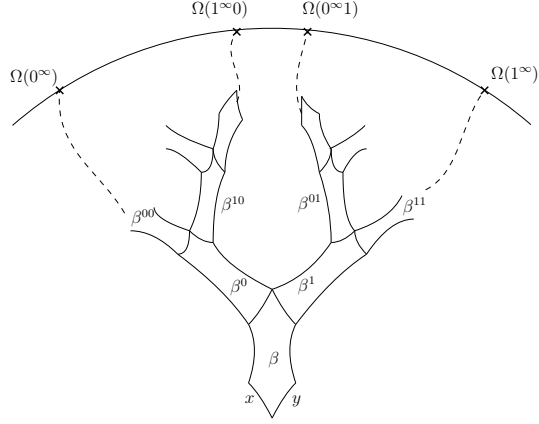}
\caption{The central interval  $C(J(x,y)) $. }
\label{fig:3}
\end{figure}

Let us consider the following interval, called {\em central} in $J(x,y)$ (see Figure~\ref{fig:3}):
\begin{equation}\label{centralInterval}
C(J(x,y)):=\big(\Omega(1^{\infty} 0;\beta(x,y)),\Omega(0^{\infty} 1;\beta(x,y))\big).
\end{equation}

The property (a) in the next result was called the {\em eventual coincidence} condition (EC) in \cite{LV} and the two properties in (b) were called the conditions (E+), (E-) when all relations in the presentation $P$ have even length.

\begin{lemma}[Eventual Coincidence of length $k$]\label{ECcentral}
If $\Theta\in C(J_1)\times C(J_2)\times\ldots\times C(J_{2N})$, then:
\begin{enumerate}
\item[(a)] For all $j\in\{1,\ldots,2N\}$,
$ \Phi_{\Theta}^{k(j)}(\theta_{j (-)}) = \Phi_{\Theta}^{k(j)}(\theta_{j (+)})$,
where $k(j)$ is the length of $\beta (x_{j-1},x_{j})$ and the notation $\Phi_{\Theta}^{k(j)}(\theta_{j(\pm)})$ stands for the iterates from the left or the right of the cutting point.
\item[(b)] If $\beta(x_{j-1},x_{j})=\{x_{j-1} x_{j_2}\cdots x_{j_k};x_{j} y_{j_2}\cdots y_{j_k} \}$ then, for each $1\leq m\leq k(j)-1$, $\Phi_{\Theta}^{m}(\theta_{j (-)})\in I_{x_{j_{m+1}}}$
and $\Phi_{\Theta}^{m}(\theta_{j (+)})\in  I_{y_{j_{m+1}}}$.
\end{enumerate}
\end{lemma}
\begin{proof}
We use again the notations $(x,y)$ for a pair of adjacent generators, as well as $\beta(x,y)$ for a minimal bigon and $k(x,y)$ for its length. Let us write
\[ \beta(x,y)=\{x\cdot x_{j_2}\cdots x_{j_k};y\cdot y_{j_2}\cdots y_{j_k}\}. \]
The extension construction defined in (\ref{extensionN}) is obtained by an infinite sequence of bigon concatenations from the initial minimal bigon $\beta(x,y)$. It is well defined by Lemma~\ref{limitExtension}, and is parameterized by $\mathbf{E}\in\{0,1\}^{\mathbb{N}}$.

At a finite step of level $n$, given by $\epsilon=(\epsilon_n,\ldots,\epsilon_1)\in\{0,1\}^{n}$, there is a bigon that we denote by $\beta^{\epsilon=(\epsilon_n,\ldots,\epsilon_1 )}$, expressed as
\begin{equation}
\beta^{\epsilon}=\{x_1^{\epsilon} x_{2}^{\epsilon}\cdots x_{k_{\epsilon}}^{\epsilon};y_1^{\epsilon} y_{2}^{\epsilon}\cdots y_{k_{\epsilon}}^{\epsilon} \}.
\end{equation}

There are two bigons at the next level, $\beta^{(0,\epsilon)}$ and $\beta^{(1,\epsilon)}$, and the concatenation operation is obtained by the identifications $x_{k_{\epsilon}}^{\epsilon}=y_{1}^{(0,\epsilon)}$ and
$y_{k_{\epsilon}}^{\epsilon}=x_{1}^{(1,\epsilon)}$.

By Proposition \ref{cylinders}, the parameter $\theta(x,y)\in\mathcal{C}_x\cap\mathcal{C}_y$ has two expressions as a geodesic ray:
\begin{equation}\label{thetaxy}
[\theta(x,y)]_{(-)} =  x\cdot x_{j_2}\cdots x_{j_k}\cdot\alpha\cdot W \mbox{ and } [\theta(x,y)]_{(+)} = y\cdot y_{j_2}\cdots y_{j_k}\cdot\alpha\cdot W,
\end{equation}
where $\alpha\in \mathcal{X}$ and $W$ is an infinite word in the alphabet $\mathcal{X}$ which might not be unique. The expression $\alpha\cdot W$ is one possible geodesic ray continuation of the two geodesics in $\beta(x,y)$. This implies in particular that $\alpha\neq (x_{j_k})^{-1}$ and $\alpha\neq (y_{j_k})^{-1}$. The two extensions of level 1 of $\beta(x,y)$ are, with the notations above:
\[ \beta^{0} = \{ x_1^{0} x_{2}^{0} \cdots x_{k_{0}}^{0};y_1^{0} y_{2}^{0}\cdots y_{k_{0}}^{0}\}\mbox{ and }\beta^{1}=\{x_1^{1} x_{2}^{1}\cdots x_{k_{1}}^{1};y_1^{1} y_{2}^{1}\cdots y_{k_{1}}^{1}\}. \]

We start by proving (b). The Definition (\ref{family}) of $\Phi_{\Theta}$ implies, by (\ref{shiftMap}), that
\[ \Phi_{\Theta} ( [\theta (x, y)]_{(-)}) =   x_{j_2} \cdots x_{j_k}\cdot\alpha\cdot W\mbox{ and }\Phi_{\Theta}([\theta(x,y)]_{(+)})=y_{j_2}\cdots y_{j_k}\cdot\alpha\cdot W. \]

By definition of the cylinder sets we obtain that $\Phi_{\Theta}([\theta(x,y)]_{(-)})\in \mathcal{C}_{x_{j_2}}$. Let us prove that
$\Phi_{\Theta}([\theta(x,y)]_{(-)})\in I_{x_{j_2}}\subset \mathcal{C}_{x_{j_2}}$. We check that
$\Phi_{\Theta}([\theta(x, y)]_{(-)})$ cannot belong to the two adjacent cylinders. The generator $x_{j_2}$ is adjacent to
$x_{j_2 \pm 1}$.

$\bullet$ On the right side, $x_{j_2+1}$ is either
$x^{-1}$ or another generator $x'$. This second possibility occurs only if $x$ and $y$ belong to a relation of length 3 (namely, $x\cdot x'\cdot y^{-1} = Id$). In this case, a direct checking shows that for all $\Theta$ we have $\Phi_{\Theta}([\theta(x,y)]_{(-)})\notin\mathcal{C}_{x^{-1}}$ or $\Phi_{\Theta}([\theta(x,y)]_{(-)})\notin\mathcal{C}_{x'}$. The same conclusion follows similarly when $x$ and $y$ belong to a relation of odd length. If $x$ and $y$ belong to a relation of even length, then the hypothesis $\theta(x,y)\in C(J(x,y))$ is needed. Indeed, from the expression of
$\theta(x,y)$ in (\ref{thetaxy}), the geodesic continuation $\alpha\cdot W$ satisfies
$\alpha\neq (x_{j_k})^{-1}$ and $\alpha\neq (y_{j_k})^{-1}$, but $\alpha$ might be adjacent to
$y_{j_k}^{-1}$ on the left, this is the worst case.
Since $\theta(x,y) \in C(J(x,y))$, the geodesic ray expression
$[\theta(x,y)]_{(+)} = y\cdot y_{j_2}\cdots y_{j_k}\cdot\alpha\cdot W$ diverges at some level from the ray
$ \Omega( 0^{\infty}\cdot1 ; \beta(x,y) )$ on its left. In this worst case, the image
$\Phi_{\Theta} ( [\theta (x, y)]_{(-)}) =   x_{j_2} \cdots x_{j_k}\cdot\alpha\cdot W$ diverges, on the left, from the ray
$ \Omega( 0^{\infty} ; \beta(x_{j_2}, x^{-1}) )$. This implies that $\Phi_{\Theta}([\theta(x,y)]_{(-)}) \notin \mathcal{C}_{x^{-1}}$.

$\bullet$ The generator on the left side of $x_{j_2}$ is $x_{j_2-1}$ which, together with $x_{j_2}$, defines a minimal bigon $\beta(x_{j_2-1},x_{j_2})$.
The point $\Phi_{\Theta}([\theta(x, y)]_{(-)})=x_{j_2} x_{j_3}\cdots x_{j_k}\cdots$ could belong to the adjacent cylinder
${\mathcal{C}}_{x_{j_2 - 1}}$ only if the initial part of $x_{j_3}\cdots x_{j_k}\cdots$ is part of the bigon
$\beta(x_{j_2-1},x_{j_2})$. This is impossible by Proposition \ref{Bigons}, since $x_{j_3}$ and the generator along
$\beta(x_{j_2-1},x_{j_2})$ following $x_{j_2}$ cannot be adjacent. Therefore, we obtain:
\begin{equation}\label{conditionb1}
\Phi_{\Theta} ( [\theta (x, y)]_{(-)}) \in I_{x_{j_2}} \textrm{ and, symmetrically, }\Phi_{\Theta} ( [\theta (x, y)]_{(+)}) \in I_{y_{j_2}}.
\end{equation}
Note that for the symmetric case we swap left and right.
The same arguments apply for all $1\leq m<k-1$, and we obtain
\begin{equation}\label{conditionbm}
\Phi_{\Theta}^m(\theta(x,y)_{(-)})\in I_{x_{j_{m+1}}} \mbox{ and  }\Phi_{\Theta}^m(\theta(x,y)_{(+)})\in I_{y_{j_{m+1}}}.
\end{equation}

A potential difficulty could occur for the iterate $m=k-1$. In this case the above arguments, on the left, using Proposition \ref{Bigons}, do not apply. At this step the extension operations in (\ref{extension0}) or (\ref{extension1}) are used. From (\ref{conditionbm}) we obtain
\begin{equation}\label{conditionbk-1}
\Phi_{\Theta}^{k-1} ( [\theta (x, y)]_{(-)}) = x_{j_k}\cdot\alpha\cdot W  \mbox{ and }\Phi_{\Theta}^{k-1}([\theta(x,y)]_{(+)})=y_{j_k}\cdot\alpha\cdot W.
\end{equation}

This implies that $\Phi_{\Theta}^{k-1} ( [\theta (x, y)]_{(-)}) \in \mathcal{C}_{x_{j_k}}$ and $\Phi_{\Theta}^{k-1} ( [\theta (x, y)]_{(+)}) \in \mathcal{C}_{y_{j_k}}$.
Observe that the generator $x_{j_k}$ is adjacent to $x_1^0$ on the left and to $x_{j_{k-1}}^{-1}$ on the right. By the argument giving (\ref{conditionb1}) we get that
\[ \Phi_{\Theta}^{k-1}([\theta(x,y)]_{(-)})\notin{\mathcal{C}}_{x_{j_{k-1}}^{-1}}. \]

The main assumption of the lemma is that $\theta(x,y)\in C(J(x,y))$. It implies (see Figure~\ref{fig:3}) that
$\Phi_{\Theta}^{k-1}([\theta(x,y)]_{(-)})\notin {\mathcal{C}}_{x_1^0}$, and thus, by definition of the intervals $I_x$, we obtain:
\begin{equation}\label{conditionbk}
\Phi_{\Theta}^{k-1} ( [\theta (x, y)]_{(-)}) \in I_{x_{j_k}} \mbox{ and }
\Phi_{\Theta}^{k-1} ( [\theta (x, y)]_{(+)}) \in I_{y_{j_k}}.
\end{equation}

This completes the argument for the point (b) of the Lemma.

The proof of (a) is direct from (\ref{conditionbk-1}) and (\ref{conditionbk}), since the definition of the map yields $\Phi_{\Theta}^{k}([\theta(x, y)]_{(-)})=
\alpha\cdot W$ and $\Phi_{\Theta}^{k}([\theta(x,y)]_{(+)})=\alpha\cdot W$, where $k=k(x,y)$.
\end{proof}

Lemma \ref{ECcentral} shows that for an open set in the parameter space, the family $\Phi_{\Theta}$ satisfies a strong dynamical property: the eventual coincidence condition of ``minimal'' length, i.e. the length of the minimal bigon. The dynamics of the maps for parameters not in this open set is a priori different and quite interesting in its own right.

\section{Centered continuations and the middle parameter}\label{middle}
In this section we define a particular cutting parameter $\Theta^0$ that we call \emph{middle}. As we will see in Section~\ref{computer},
this choice of $\Theta^0$ allows us to perform particularly simple computations in order to get the volume entropy of the presentation. Many other parameters would allow the
same computations, in particular many parameters in the \emph{central} interval defined in Section~\ref{central}, but $\Theta^0$ is very simple to define.

For each $1\le i\le 2N$, we define the \emph{opposite} to $x_i$ as the symbol $x_{i+N}$ (recall that the subindexes in
the generating set $\{x_1,x_2,\ldots,x_{2N}\}$ are taken modulo $2N$, with $1,\ldots,2N$ as representatives of
the classes modulo $2N$).

\begin{figure}
\centering
\includegraphics[scale=0.5]{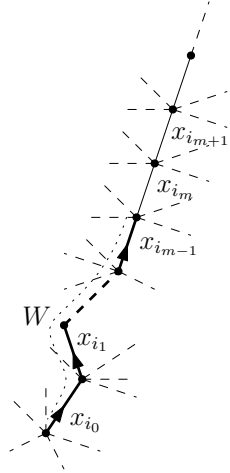}
\caption[fck]{The centered continuation of a geodesic $W$}\label{centered}
\end{figure}

Let $g\in G$ and let $W=x_{i_0}x_{i_1}\cdots x_{i_{m-1}}$ be a geodesic expression of $g$. The infinite word $Wx_{i_m}x_{i_{m+1}}\cdots$ such
that $x_{i_j}$ is the opposite to $(x_{i_{j-1}})^{-1}$ for each $j\ge m$ will be called the \emph{centered continuation of} $W$. Note that,
from the vertex $g$ of the Cayley graph, the last generator $x_{i_{m-1}}$ in the word $W$ reads as the edge labeled $(x_{i_{m-1}})^{-1}$.
See Figure~\ref{centered}. From Lemma~2.9 in \cite{L}, the centered continuation of a geodesic expression of $g$ is a geodesic ray.
By definition, the corresponding point in $S^1$ belongs to $\mathcal{C}_g$. It will be called the \emph{$W$-middle point}.

The next result will be used to prove one of the two inequalities stated in Proposition~\ref{SigmaVsX}, a
central result of Section~\ref{inequalities}.

\begin{lemma}\label{continuation}
Let $g\in G$ and let $W=x_{i_0}x_{i_1}\cdots x_{i_{m-1}}$ be a geodesic expression of $g$. Let $Wx_{i_m}x_{i_{m+1}}\cdots$ be the centered continuation of $W$
and let $\zeta \in \partial G$ be the $W$-middle point. Then, any geodesic ray representing $\zeta$ has the form $\hat{W}x_{i_m}x_{i_{m+1}}\cdots$, where $\hat{W}$ is a
geodesic expression of $g$.
\end{lemma}
\begin{proof}
Set $L=Wx_{i_m}x_{i_{m+1}}\cdots$ and, for $0\le i\le m$, let us denote by $v_i$ the vertex at distance $i$ from $\Id$ inside the geodesic
segment $W$. So, $v_0=\Id$ and $v_m=g$. Let us assume, by way of contradiction, that there exists a geodesic ray $R$ whose vertex
at distance $m$ from the identity is $g'\ne g$. Assume without loss of generality, just to fix ideas and fit to the notation $L/R$,
that $g'$ is on the right of $g$ in the obvious orientation induced in the Cayley graph by the clockwise orientation
of the circle. Let $j$ be the largest index in $\{0,1,\ldots,m\}$ such that $v_j\in R$. We are assuming that $j<m$. Now we have two cases.

Assume first that $L$ and $R$ intersect beyond $v_j$. Let $w$ be the vertex in $L\cap R$ at least distance from $\Id$ larger than $m$.
Then, the segments $\bar{L}\subset L$ and $\bar{R}\subset R$ connecting $v_j$ and $w$ are geodesics and $\bar{L}\cap\bar{R}=\{v_j,w\}$.
By Lemma~\ref{jer2}, they are the left and right sides of a staircase between $v_j$ and $w$. Remarks~\ref{restriction} and
\ref{no3} imply then that, given any vertex in $\bar{L}$, there can be at most one edge of the Cayley graph
contained in the bounded component $K$ of $\mathbb{R}^2\setminus(\bar{L}\cup\bar{R})$. But, by definition of the centered
continuation, $x_{i_m}$ is the edge opposite to $(x_{i_{m-1}})^{-1}$. Since $2N\ge6$, this implies that there are at least two
edges departing from $g\in\bar{L}$ and contained in $K$ (see Figure~\ref{centered}), a contradiction that proves the result in this case.

Assume now that $L$ and $R$ does not intersect beyond $v_j$. Therefore, the infinite rays $\bar{L}\subset L$ and $\bar{R}\subset R$
starting at $v_j$ define a bigon ray based at $v_j$. Then, from Lemma~\ref{jer4} it follows that $\bar{L}$ and $\bar{R}$ are
respectively the left and right sides of an infinite extension of bigons based at $v_j$, and we get a contradiction using
exactly the same arguments as in the previous case.
\end{proof}

The centered continuation allows us to define a particular choice $\Theta$ of the cutting parameter so that the dynamics of
any cutting point by $\Phi_{\Theta}$ is specially simple.

For each $1\le i\le 2N$, let $L_i^l$ and $L_i^r$ be the left and
right geodesic segments of length $k(i):=k(x_{i-1},x_i)$ defining the
minimal bigon $\beta(x_{i-1},x_i)$. The cutting parameter $\Theta^0 =(\theta^0_1,\theta^0_2,\ldots,\theta^0_{2N})$ such that
$\theta^0_i$ is the $L_i^l$-middle point for $i\in\{1,2,\ldots,2N\}$ will be called the \emph{middle parameter}.

\begin{remark}\label{ECcentralInterval}
Let $\Theta^0=(\theta^0_1,\theta^0_2,\ldots,\theta^0_{2N})$ be the middle parameter. Then, for each $1\le j\le 2N$,
the cutting point $\theta^0_j$ belongs to the central interval $C(J_j)$ defined in Section~\ref{central}. Indeed, if $z$
is the last generator in the word $L_j^l$, then the next generator along the centered continuation is the opposite to $z^{-1}$,
that does not belong to the bigon extension (\ref{Extension}) of $\beta(x_{j-1},x_j)$.
\end{remark}

In view of Remark~\ref{ECcentralInterval}, the first $k(j)$ iterates from the left and the right of $\theta^0_j$ by
$\Phi_{\Theta^0}$ are controlled by properties (a) and (b) of Lemma~\ref{ECcentral}. The subsequent iterates satisfy the
property stated in the next result.

\begin{lemma}\label{ECcentralInterval2}
Let $\Theta^0=(\theta^0_1,\theta^0_2,\ldots,\theta^0_{2N})$ be the middle parameter. Then, for each $1\le j\le 2N$,
$\Phi_{\Theta^0}^m(\theta^0_j)$ is not a cutting point for any $m\ge k(j)$.
\end{lemma}
\begin{proof}
The definition of the centered continuation implies that $\Phi_{\Theta^0}^{k(j)}(\theta^0_j)$ belongs to an interval
$I_{x_{n}}$ and more precisely to $I_{x_{n}}\setminus (J_{x_{n -1}}\cup J_{x_{n + 1}})$.
Indeed, along the geodesic ray $[\theta^0_j]$ defining $\theta^0_j$, after the initial bigon $\beta(x_i,x_{i+1})$, there is never
two adjacent generators and thus there is never a subword of a bigon. Therefore, all iterates $\Phi_{\Theta^0}^{m} (\theta^0_j)$,
for $m \ge k(j)$, belong to some interval $I_{x_{r(m)}}\setminus (J_{x_{r(m) -1}} \cup J_{x_{r(m)+ 1}})$ for some $r(m) \in \{1, 2, ..., 2N\}$.
In consequence, $\Phi_{\Theta^0}^{m}(\theta^0_j)$ is not a cutting point for all $m \ge k(j)$.
\end{proof}

\section{Comparison of entropy functions}\label{inequalities}
We start this section by recalling how the topological entropy of a piecewise monotone discontinuous map of the circle can be computed as
the exponential growth rate of a function $X_m$ that counts the number of itineraries of length $m$ associated to a given partition
of $S^1$ such that the map is monotone on each interval of the partition (Lemma~\ref{topEntropy}). Then we prove Proposition~\ref{SigmaVsX},
that compares the dynamical counting function $X_m$ with the growth function $\sigma_m$ defined in (\ref{sigma-n}). As a corollary
we get Theorem~\ref{entropies}, that states the equality between the topological entropy of any Bowen-Series-like map and the volume entropy of
the group presentation.

The maps in the family $\Phi_{\Theta}$ are defined by a partition
$S^1=\bigcup_{j=1}^{2N} I_j$, where the intervals
$I_j=[\theta_j,\theta_{j+1})$ depend on the cutting parameter and the
restriction of $\Phi_{\Theta}$ to each $I_j$ is a homeomorphism. The intervals
$I_j$ will be called \emph{basic intervals} from now on. The
generating set $\mathcal{X}$ of the group $G$ is in bijection with
the set of basic intervals, for all parameters $\Theta$.
Each point $z\in S^1$ has an orbit which is described by an
itinerary, i.e. an ordered sequence
\[ It(z):=(j_0,j_1,\ldots,j_m,\ldots) \mbox{ such that } \Phi_{\Theta}^m(z) \in I_{j_m} \forall\, m\ge0. \]
For $m\in\mathbb{N}$, the \emph{itinerary intervals of level $m$} are defined as
\[ I_{j_0,j_1,\ldots,j_{m-1}}:=I_{j_0}\cap\Phi^{-1}_{\Theta}(I_{j_1})\cap\ldots\cap\Phi^{-(m-1)}_{\Theta}(I_{j_{m-1}}) \]
and the dynamical counting function is given by:
\begin{equation*}
X_{m} \textrm{ is the number of non-empty intervals } I_{j_0,j_1,\ldots,j_{m-1}}\textrm{ of level } m.
\end{equation*}

Since each $\Phi_{\Theta}$ is a piecewise homeomorphism, there is at most one
connected component in $I_{j_0}$ such that $\Phi_{\Theta}(I_{j_0})\cap I_{j_1}\neq\emptyset$.
In consequence, the itinerary intervals are well defined.

\begin{remark}\label{new}
Note that if $I_{j_0,j_1,\ldots,j_{m-1}}$ is a non-empty interval then there exists $z\in S^1$ admitting a $\Phi_{\Theta}$-orbit that visits
$I_{x_{j_0}},I_{x_{j_1}},\ldots,I_{x_{j_{m-1}}}$. Therefore, if $I_{j_0,j_1,\ldots,j_{m-1}}\ne\emptyset$, then the path
$x_{j_0}x_{j_1}\cdots x_{j_{m-1}}$ is a geodesic path in the Cayley graph and the element $g=x_{j_0}x_{j_1}\ldots x_{j_{m-1}}$ has length $m$.
\end{remark}

The itinerary intervals satisfy the following properties, that follow immediately from the definitions:
\begin{enumerate}
\item[(a)] $\forall z\in I_{j_0,j_1,\ldots,j_m},\, \Phi^{r}_{\Theta}(z)\in I_{j_r}$ for $0\le r\le m$.
\item[(b)] $\forall z \in  I_{j_0,j_1,\ldots,j_m},\, \Phi^{m+1}_{\Theta}(z)=g^{-1}z$ where $g=x_{j_0}x_{j_1}\ldots x_{j_m}$.
\item[(c)] $I_{j_0,j_1,\ldots,j_m} $ is a subinterval of $I_{j_0,j_1,\ldots,j_{m-1}}$.
\end{enumerate}

\begin{lemma}\label{topEntropy}
With the above notations, for each cutting parameter $\Theta\in J_1 \times \ldots \times J_{2N}$ the following equality holds: $h_{\textrm{top}}(\Phi_{\Theta})=\lim_{m \rightarrow\infty}\frac{1}{m}\log(X_m)$.
\end{lemma}
\begin{proof}
The original definition of the topological entropy \cite{AKM} uses a counting function for coverings of the
space. The definition holds for continuous self-maps on compact spaces. Bowen in \cite{Bo1} extended this
notion to non necessarily continuous self-maps. In \cite{MZ} Misiurewicz and Ziemian showed
that, for piecewise continuous and piecewise monotone interval maps, the entropy can be computed from \emph{mono-covers},
namely partitions of the interval by subintervals such that the map restricted to each element in the partition is continuous
and monotone. With this partition, we denote by $Y_m$ the number of different itineraries of length $m$. Then,
Misiurewicz and Ziemian showed that the entropy of the map coincides with
\[ \lim_{m \rightarrow\infty} \frac{\log Y_m}{m}.\]

It is straightforward to adapt the previous result to the case of piecewise continuous and piecewise monotone maps of the circle.
\end{proof}

In view of Lemma~\ref{topEntropy} and the definition (\ref{volentropy}) of the volume entropy, that uses the counting function
$\sigma_m$ introduced in (\ref{sigma-n}), we need to compare the rates of increasing with $m$ of $X_m$ and $\sigma_m$. This
is the aim of the next result.

\begin{proposition}\label{SigmaVsX}
With the above notations, the following inequalities are satisfied for each parameter $\Theta\in J_1\times J_2\times\ldots\times J_{2N}$:
\[ \sigma_m \le X_m\le m\sigma_m. \]
\end{proposition}
\begin{proof}
Let $g\in G$ be an element of length $m$. We will show that the number of nonempty intervals $I_{j_0,\ldots,j_{m-1}}$ such that
$x_{j_0}x_{j_1}\ldots x_{j_{m-1}}=g$ is between 1 and $m$.

Consider a geodesic word $W=x_{i_0}x_{i_1}\cdots x_{i_{m-1}}$ representing $g$. Let $Wx_{i_m}x_{i_{m+1}}\cdots$ be the
centered continuation of $W$ and let $\rho\in S^1$ be the corresponding $W$-middle point. For $0\le i<m$,
let $I_{x_{j_i}}$ be the only interval of the partition $\{I_{x_j}\}_{j=1}^{2N}$ of $S^1$ containing $\Phi_{\Theta}^i(\rho)$.
Then, the itinerary interval $I_{j_0,j_1,\ldots,j_{m-1}}$ is nonempty and, by Remark~\ref{new}, $x_{j_0}x_{j_1}\cdots x_{j_{m-1}}$
is a geodesic path in the Cayley graph. Then, from Lemma~\ref{continuation} it follows that $x_{j_0}x_{j_1}\ldots x_{j_{m-1}}=g$.
This shows that $\sigma_m \le X_m$.

\begin{figure}
\centering
\includegraphics[scale=0.6]{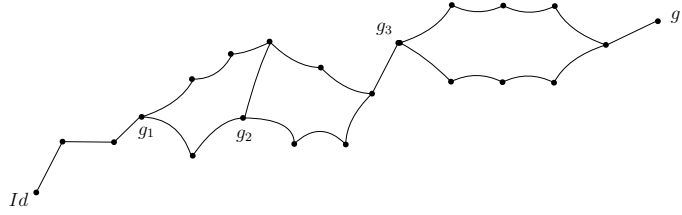}
\caption[fck]{The points $g_i$ as defined in the proof of Proposition~\ref{SigmaVsX}.}\label{gallery2}
\end{figure}

Let us prove now the second inequality. Clearly, if $I_{j_0,\ldots,j_{m-1}}\ne\emptyset$
then the path in the Cayley graph corresponding to the word $g:=x_{j_0}x_{j_1}\cdots x_{j_{m-1}}$
is a geodesic between the vertices $\Id$ and $g$. So, it is
contained in the set of all geodesics joining $\Id$ and $g$ which, by
Lemma~\ref{totes}, is a gallery that we will denote by
$H_{\Id,g}$. From the definition of a gallery, there are at most two
edges in $H_{\Id,g}$ departing from the identity. This implies that the cylinder set
${\mathcal{C}}_g$ defined in (\ref{cylinderx}) of all infinite geodesic rays
passing through $g$ is contained in a set $K$ that is either
a basic interval or the union of two adjacent basic intervals.
Since $H_{\Id,g}$ is a gallery and the length of $g$ is $m$, $H_{\Id,g}$ contains
$k$ minimal bigons for some $0\le k<m$. We denote by $g_1,g_2,\ldots,g_k$ the vertices
where such bigons are based at. See Figure~\ref{gallery2}. Let $\theta_{i_1},\theta_{i_2},\ldots,\theta_{i_k}$ be
the corresponding cutting points associated to each of these bigons.
For $1\le j\le k$, set $\tilde\theta_j:=g_j\theta_{i_j}$. Note that the cardinality
of the set $\{\tilde{\theta}_j\}_{j=1}^k$, that is contained in $K$ by construction, is between 1 and $k$
(see Figure~\ref{gallery2} for an example). Therefore,
\[ K\setminus \{\tilde \theta_1,\ldots,\tilde\theta_k\}\]
is the union of at most $k+1$ disjoint intervals that we denote
by $K_1,K_2,\ldots,K_{k+1}$.

We claim that in each of such intervals there is at most one nonempty interval $I_{j_0,\ldots,j_{m-1}}$
satisfying $x_{j_0}x_{j_1}\cdots x_{j_{m-1}} = g$. Indeed, assume on the contrary that for some $l$ there are two
nonempty intervals $I_{j_0,\ldots,j_{m-1}}$ and $I_{h_0,\ldots,h_{m-1}}$ contained in $K_l$ such that
$x_{j_0}x_{j_1}\cdots x_{j_{m-1}} = x_{h_0}x_{h_1}\cdots x_{h_{m-1}} = g$.

The corresponding geodesic segments are thus contained in the gallery $H_{Id,g}$. Let $s$ be the least index such that
$j_s\ne h_s$. Then, $x_{j_0}x_{j_1}\cdots x_{j_{s-1}} = x_{h_0}x_{h_1}\cdots x_{h_{s-1}} = g_r$ for some $g_r$ being the
origin of a minimal bigon $\beta_{g_r}\subset H_{\Id,g}$ associated to some cutting point $\theta_{i_r}$.
The bigon $\beta_{g_r}$ starts at $g_r$ by the edges labeled $x_{j_s}$ and $x_{h_s}$. Therefore the two itinerary intervals
$I_{j_0,\ldots,j_{m-1}}$ and $I_{h_0,\ldots,h_{m-1}}$ contained in $K_l$ have a preimage of order $s$ of the cutting point
$\theta_{i_r}$ between them, i.e. $\Phi_{\Theta}^{-s}(\theta_{i_r}) = g_r (\theta_{i_r})$. Thus
$g_r (\theta_{i_r}) = \tilde{ \theta}_{i_r}$ is contained in the interior of $K_l$, a contradiction that proves the claim.

Since $k+1\le m$, it follows that there are at most $m$ nonempty intervals $I_{j_0,\ldots,j_{m-1}}$ such that
$x_{j_0} x_{j_1}\ldots x_{j_{m-1}} = g$.
\end{proof}

Finally, we get the desired equality between the entropies.

\begin{theorem}\label{entropies}
Let $\Sigma$ be a compact closed hyperbolic surface with a geometric presentation $P$ for the group $G:=\pi_1(\Sigma)$. Then,
\[ h_{vol}(G,P)=h_{top}(\Phi_{\Theta})\mbox{ for all }\Theta\in J_1 \times\ldots\times J_{2N}. \]
\end{theorem}
\begin{proof}
It follows directly from (\ref{volentropy}), Lemma~\ref{topEntropy} and Proposition~\ref{SigmaVsX}.
\end{proof}

\section{Proof of Theorem~\ref{mainTh}}\label{milnor}
The proof of Theorem~\ref{mainTh} strongly uses the Milnor-Thurston theory \cite{MT,AM2} of kneading invariants. We start this
Section by reviewing the main ingredients of the theory. Then we see how the results of Milnor and Thurston apply to our context. Finally we prove Theorem~\ref{mainTh}.

Kneading theory of Milnor and Thurston was stated for piecewise monotone
and continuous maps of the interval. For such a map
$\map{f}{I}$, the lap number $l(f)$ is the number
of maximal intervals of monotonicity of $f$. Each point separating
different maximal intervals of monotonicity is called a \emph{turning
point}. In this situation, the limit $s(f):=\lim (l(f^n))^{\frac{1}{n}}$ exists
and is called the \emph{growth number} of $f$. It is well known \cite{ms} that the
topological entropy of $f$ is precisely $h_{top}(f)=\log(s(f))$.

Let $I_1,I_2,\ldots,I_l$ be the laps associated to the map $f$, numbered
in the natural ordering along $I$, and let $v_1<v_2<\ldots<v_{l-1}$ be the corresponding
turning points. The \emph{address} $A(x)$  for a point $x\in I$  will be the
formal symbol $I_j$ if $x\in I_j$ or the formal symbol $C_j$ if
$x=v_j$. The \emph{itinerary} of $x$ is the infinite sequence of symbols
\[ I(x)=A(x),A(f(x)),A(f^2(x)),\ldots \]

To each address $A(x)$ we assign the number
\[ \epsilon(A(x)):=
\left\{
\begin{array}{ll}
1 & \hbox{if $f$ is increasing at $x$} \\
-1 & \hbox{if $f$ is decreasing at $x$} \\
0 & \hbox{if $x$ is a turning point}
\end{array}
\right.
\]
and define a new infinite sequence $\omega(x)=(\omega_0(x),\omega_1(x),\ldots,\omega_n(x),\ldots)$ given by
\[ \omega_j(x)=\epsilon_j(x)A(f^j(x)), \mbox{ where }
\epsilon_j(x)=\epsilon(A(x))\epsilon(A(f(x)))\ldots\epsilon(A(f^{j-1}(x))). \]

Note that
\[\epsilon_j(x)= \left\{
\begin{array}{ll}
1 & \hbox{if $f^{j-1}$ is increasing at $x$} \\
-1 & \hbox{if $f^{j-1}$ is decreasing at $x$} \\
0 & \hbox{if $f^k(x)$ is a turning point for some $k<j-1$}.
\end{array}
\right.\]

For any $x\in I$ and $k\in\mathbb{N}$ it is clear that $\omega_k(y)$ is constant
in a neighborhood $(x-\delta_k,x)$. We denote the corresponding symbol by $\omega_k(x^-)$. In
the symmetric way we define $\omega_k(x^+)$. Finally, we define the
\emph{jump at} $x$ as the series
\[\Omega_x(t)=\sum_{i=0}^{\infty}\left(\omega_i(x^+)-\omega_i(x^-)\right)t^i.\]
Note that if for all $n\in N,\,\,f^n(x)$ is not a turning point, then $\Omega_x(t)$ is identically zero.
The \emph{$j$-kneading invariant} is defined as the jump at the $j$-th turning point, that is
\[ \nu_j(t)=\Omega_{v_j}(t)=\sum_{i=0}^{\infty}\left(\omega_i(v_j^+)-\omega_i(v_j^-)\right)t^i. \]
It is clear that $\nu_j(t)=N_{j 1}I_1+\ldots N_{j l}I_l$, for any $j\in \{1,2,\ldots,l-1\}$, where
$N_{j i}$ belongs the ring of formal power series with integer coefficients $\mathbb{Z}[[t]]$.

At this moment we define the $(l-1)\times l$ kneading matrix $M_f$
whose entries are the $N_{ij}$ power series. Let us denote by
$D_i$ the determinant of the
matrix obtained by deleting the $i$-th column of the kneading matrix. The
following is the first result of the above construction.

\begin{theorem}\label{kneading} With the above notation, the series $D(t)=(-1)^{i+1}\frac{D_i}{1-\epsilon(I_i)t}$
is a fixed element of $\mathbb{Z}[[t]],$ independently of the choice of $i$.
It has leading coefficient 1, converges for $|t|<1$ and $D(t)\ne 0$
if $|t|<1/s(f)$. Moreover, the first zero of $D(t)$ in the interval
$[0,1)$ occurs at $\frac{1}{s(f)}$ if $s(f)\ge 1$, while $D(t)$ has no
zeros in the open unit disk if $s(f)=1$.
\end{theorem}

In short, Theorem~\ref{kneading} says that the entropy of a map can
be computed knowing the itineraries of the turning points.

Another important result from the Milnor-Thurston theory is the
semi-conjugacy of any piecewise monotone map $f$ of the interval to a
piecewise affine map with constant absolute value of the slopes.
Next we give a brief explanation on how to construct the semi-conjugacy.
See \cite{MT} for details. First of all consider the series
\[ L(t)=\sum_{i=0}^\infty l(f^n)t^n. \]
Clearly $L(t)$ converges for $t<1/s(f)$ and has a pole at $t=1/s(f)$. For any subinterval $J$ of $I$,
we define \[ L[J](t)=\sum_{i=0}^\infty l(f^n\vert_J)t^n, \]
where $l(f^n\vert_J)$ denotes the number of laps of $f^n$ that intersect
$J$. Now we define \[ \mu(J):=\lim_{t\to 1/s(f)}\frac{L[J](t)}{L(t)}.\]
It is not difficult to see that $\mu$ satisfies the following properties:
\begin{enumerate}
\item[(i)] $0\le \mu(J)\le 1$
\item [(ii)] If $J_1\cap J_2$ reduces to a single point, then $\mu(J_1\cup J_2)=\mu(J_1)+\mu(J_2)$
\item [(iii)] The number $\mu(J)$ depends continuously on the endpoints of $J$
\item [(iv)] If $J$ is contained in single lap of $f$, then $\mu(f(J))=s(f)\mu(J)$
\end{enumerate}

Properties (i--iv) imply that $\mu$ gives rise to a measure on the interval $I$. Now consider the map $\map{h}{I}[[0,1]]$
defined as
\begin{equation}\label{defh}
h(x)=\mu([a,x]), \tag{$\star$}
\end{equation}
where $a$ is the left endpoint of $I$. From properties (i--iv) above it follows that $h$ is a continuous and
non-decreasing map from $I$ onto $[0,1]$. Moreover, it is not difficult to see that if
$f$ is topologically transitive then $h$ is a homeomorphism (see for instance Proposition~4.6.9 of \cite{alm}).
Lastly we have the following result, which is essentially Theorem~7.4 of \cite{MT}.

\begin{theorem}\label{conju}
Let $\map{f}{I}$ be a piecewise monotone continuous map of an interval $I$ and
let $\map{h}{I}[[0,1]]$ be the corresponding non-decreasing map defined by (\ref{defh}).
Then, there exists a unique piecewise linear map $\map{F}{[0,1]}$ with constant slope $\pm s(f)$ such that
\[ F(h(x))=h(f(x)). \]
Moreover, if $f$ is topologically transitive, then $h$ is a homeomorphism.
\end{theorem}

In \cite{AM1}, Theorems~\ref{kneading} and \ref{conju} have been generalized to piecewise monotone maps $f$ of
degree one of the circle by considering its lifting $\tilde{f}$ and the map $\map{\widehat{f}}{[0,1]}$ defined by
$\widehat{f}(x)=\tilde{f}(x)-E(\tilde{f}(x))$, where $E(y)$ denotes the integer part of $y$. The map $\widehat{f}$
is also piecewise monotone but it is discontinuous. Therefore, it is necessary to consider the discontinuity points
as turning points. In \cite{AM2}, the same results have been also generalized to discontinuous piecewise
monotone maps of the circle. Now all the discontinuities are considered as turning points.

Let us go back to our Bowen-Series-like maps.
The next result states that, when $\Theta$ is the \emph{middle parameter} as defined in Section~\ref{middle},
the jump series of $\Phi_{\Theta}$ in all turning points are in fact polynomials.

\begin{lemma}\label{essencial} Let $\Theta^0=(\theta^0_1,\theta^0_2,\ldots,\theta^0_{2N})$ be the middle parameter.
Then, the jump series $\Omega_{\theta^0_i}(t)$ for the map $\Phi_{\Theta^0}$ is a polynomial of degree at most
$k(x_{i-1},x_i)-1$ for every $1\le i\le 2N$. In consequence, the kneading invariants of the map $\Phi_{\Theta^0}$ are
integer polynomials.
\end{lemma}
\begin{proof}
We denote by $\tilde\Phi_{\Theta^0}$ a lifting of
$\Phi_{\Theta^0}$. Identifying our topological $S^1$ with the unit circle in the complex plane,
let $z_1,z_2,\ldots,z_{2N}\in [0,1)$ be such
that $e^{2\pi z_i}=\theta^0_i$. Without loss of generality we can
assume that $z_1=0$. For $i=1,2,\ldots,2N$, we set
$I_i:=(z_i,z_{i+1})$ (where $z_{2N+1}=0$). As explained, we will
consider the map $\widehat{\Phi_{\Theta^0}}$, denoted by $\Phi$ from now on
in order to simplify the notation. Observe that some of the intervals
$I_i$ may split into two laps denoted by $I_i^l$ and $I_i^r$ (standing
for \emph{left} and \emph{right}) separated by
one pre-image of $0$, denoted by $m_i$.
The map $\Phi$ has, say, $l$ laps and $l-1$ turning points,
namely $z_2,\ldots,z_{2N}$ plus the set of points $m_i$. See Figure~\ref{km2}
in Section~\ref{computer} for a particular example, where the crosses in the horizontal
axis mark the positions of the points $m_i$.

Recall now the notation $\mathcal{X} = \{x_1,x_2,\ldots,x_{2N}\}$ introduced in Definition~\ref{standing} for
the elements of the generating set. For each $1\le j\le 2N$, let $\beta_j$ be the minimal bigon
$\beta(x_{j-1},x_j)$ and let us denote its length by $k_j$. Each bigon
$\beta_j$ is defined by two geodesics $L_j^l=x_{j-1}x_{j_2}\cdots x_{j_{k_j}}$ and $L_j^r=x_jy_{j_2}\cdots y_{j_{k_j}}$.
By Remark~\ref{ECcentralInterval}, each cutting point $\theta^0_j$ belongs to the central interval $C(J_j)$. In
consequence, we can apply Lemma~\ref{ECcentral}(a). It follows that
\[ \lim_{h\to 0^+} \Phi^{k_j}(z_j+h)=\lim_{h\to 0^-} \Phi^{k_j} (z_j+h)=:\pi_j. \]

Now we claim that $\epsilon_{k_j}(z_j^+)=\epsilon_{k_j}(z_j^-)$. Indeed, since $\theta^0_j \in C(J_j)$ it can be written as two different geodesic rays $L_j^lW,L_j^rW$.
At a symbolic level, $\Phi^{k_j}(\theta_j^0)_{(+)}$ corresponds to shifting the first $k_j$ letters of $L_j^lW$, so giving $W$
(see Lemma \ref{ECcentral}).
Now, if we apply $(\Phi^{-1})^{k_j}_{(-)}$,
which corresponds to the left multiplication by $L_j^r$, we get $L_j^rW$. So, the composition $(\Phi^{-1})^{k_j}_{(-)}\circ\Phi^{k_j}_{(+)}$ is the identity map on a neighborhood of $\theta^0_j \in C(J_j)$.
At a group theoretic level, this identity is nothing but a relation in the group, seen as a subgroup of $\textrm{Homeo}(S^1)$.
In consequence, $\Phi^{k_j}_{(-)}$ and $\Phi^{k_j}_{(+)}$ are either both increasing or both decreasing. From this fact and the definition of the terms of the jump series,
the claim follows.

From the previous claim we get that
\[ \Omega_{z_j}(t)=\sum_{i=0}^{k_{j-1}}\left(\omega_i(z_j^+)-\omega_i(z_j^-)\right)t^i+\epsilon_{k_j}(z_j)t^{k_j}\Omega(\pi_j).\]
Since the orbit of $\pi_j$ does not visit any cutting point by Lemma \ref{ECcentralInterval2} then $\Omega (\pi_j)=0$.
\end{proof}

With all these ingredients we are ready to prove Theorem~\ref{mainTh}.

\begin{proof}[Proof of Theorem~\ref{mainTh}]
The equality between the entropies stated in (a) follows from Theorem~\ref{entropies}. The fact that $1/\lambda$ is the smallest root in $(0,1)$ of an integer polynomial follows
from Lemma~\ref{essencial} and Theorem~\ref{kneading}, which states that for the middle parameter the entries of the kneading matrix are integer polynomials.

The proof of statement (b) follows from Proposition~\ref{transitive} and the generalization of Theorem~\ref{conju} to circle maps \cite{AM2,AM1}.
\end{proof}

\section{Effective computations and examples}\label{computer}
In this section we give some examples of computation of the volume entropy associated to a geometric presentation. The whole process is based on Theorem~\ref{mainTh} and the Milnor-Thurston theory outlined in Section~\ref{milnor}. Moreover, the power of such techniques has allowed us to organize the computation in a systematic way to produce an algorithm that takes as input any presentation $P$ given by a list $R$ of relations and prints out the associated volume entropy provided that $P$ is geometric (otherwise, the program reports that $P$ is not geometric). The procedure has been written in Maple and Maxima languages and is freely available to the scientific community upon request to the authors. Next we show the source code of the main program, just because it clearly sets out the main steps of the computation.

{\scriptsize \begin{verbatim}

VolumeEntropy := proc(R)
    # Computes the volume entropy of the presentation R if R is geometric
    # or reports that R is not geometric otherwise
    # USES: LinearAlgebra library
    # USES: CheckRelations,CyclicOrdering,MinimalBigons,KneadingMatrix procedures
    local co,bigons,i,A:
    if CheckRelations(R)=false then
      print(`The relations do not satisfy the syntax conventions`): return:
    end if:
    co:=CyclicOrdering(R):
    if co=false then print(`The presentation is not geometric`): return: end if:
    bigons:=MinimalBigons(R,co):
    A:=KneadingMatrix(co,bigons,R):
    i:=Determinant(DeleteColumn(A,1)):
    return 1/min(fsolve(i,t=0..1));
  end proc:

\end{verbatim}}

See Table~\ref{outputs} for some examples of execution of the program. In the third column, it is shown the polynomial factor of the kneading determinant with roots in $[0,1)$.

\begingroup
\begin{table}
\centering
\setlength{\tabcolsep}{2pt}
\renewcommand*{\arraystretch}{1.3}
\begin{tabular}{|c|c|c|} \hline
{\small Presentation (relations)} & {\small Program output} & {\small Polynomial} \\
\hline
{\small $[ acde\bar{d}\bar{b},\bar{e}\bar{c}b\bar{a} ]$} & {\small $\log(8.50591006)$} & {\small $t^4-7t^3-12t^2-7t+1$} \\
\hline
{\small $[ acde\bar{b},\bar{d}e\bar{c}b\bar{a} ]$} & {\small $\log(8.78515105)$} & {\small $t^4-8t^3-6t^2-8t+1$} \\
\hline
& & {\small $t^{20}-4t^{19}-44t^{18}-122t^{17}$} \\
& & {\small $-206t^{16}-280t^{15}-381t^{14}-484t^{13}$} \\
{\small $[ aba\bar{c}d,ce^2,dbf^2 ]$} & {\small $\log(9.91984307)$} & {\small $-579t^{12}-606t^{11}-606t^{10}-606t^9$} \\
& & {\small $-579t^8-484t^7-381t^6-280t^5$} \\
& & {\small $-206t^4-122t^3-44t^2-4t+1$} \\
\hline
{\small $[ aihlk\bar{c}a,\bar{c}e^2,$} & & \\
{\small $dbf^2k,g\bar{h}j^2,idgb\bar{l} ]$} & {\small Non geometric} & -- \\
\hline
&  & {\small $t^{20}-13t^{19}-80t^{18}-149t^{17}$} \\
{\small $[ aia\bar{c}h,ce^2,dbf^2,$} & & {\small $-187t^{16}-196t^{15}-252t^{14}-348t^{13}$} \\
{\small $g\bar{h}j^2,idgb ]$} & {\small $\log(17.9527833)$} & {\small $-370t^{12}-426t^{11}-312t^{10}-426t^9$} \\
& & {\small $-370t^8-348t^7-252t^6-196t^5$} \\
& & {\small $-187t^4-149t^3-80t^2-13t+1$} \\
\hline
\end{tabular}
\caption[fck]{Some outputs of the algorithm.}\label{outputs}
\end{table}
\endgroup

\subsection{Checking relations}
The procedure {\bf CheckRelations(R)} checks whether $R$ verifies the following syntax conventions. The set $R$ of relations must be a list of lists of integers. Any generator is entered as a positive integer $k$, while its inverse is $-k$. The set of absolute values of all integers in $R$ has to be $\{1,2,\ldots,n\}$ for some $n\ge3$. By Remark~\ref{restriction}, we exclude presentations having relations of length 2 and also the presentation having 3 generators and 2 relations of length 3. As an example, next we show the inputs satisfying the conventions of the program for the following presentations (from now on, the inverse of an element $a$ will be written as $\bar{a}$):
\[ P_1 = \langle a,b,c,d\, |\,adac,cbdb\rangle\longrightarrow R=[ [1,4,1,3],[3,2,4,2] ] \]
\[ P_2 = \langle a,b,c,d,e\, |\,abc,ce\bar{a},b\bar{c}d^2\rangle\longrightarrow R=[ [1,2,3],[3,5,-1],[2,-3,4,4] ] \]
\[ P_3 = \langle a,b,c,d\, |,aba\bar{b}d,c^2d\rangle\longrightarrow R=[ [1,2,1,-2,4],[3,3,4] ] \]

\subsection{Computing the cyclic ordering}
Once we have a syntactically admissible list of relations $R$, we have to test whether it corresponds to a geometric presentation. This task is assigned to the procedure {\bf CyclicOrdering(R)}, based on properties (a,c,d) of Lemma~\ref{co}, which is in fact a characterization of the geometric presentations of co-compact hyperbolic surface groups. As a first test of geometricity, we check whether each generator appears twice (with $+$ or $-$ sign) in $R$ (Lemma~\ref{co}(c)). This test excludes the presentation $P_2$ above as non geometric, while $P_1$ and $P_3$ pass the test. It turns out, however, that $P_3$ is geometric while $P_1$ is not. The reason is that $P_3$, unlike $P_1$, satisfies property (a): it is possible to find a cyclic ordering in the generating set that is preserved by the action of the group. Let us briefly sketch how to decide whether such a cyclic ordering exists for a presentation $P=\langle X|R\rangle$ and, in the affirmative, how to find it.

First, we fix the identity element $\Id$ of the group as the base vertex of the Cayley 2-complex $\textrm{Cay}^2(G,P)$. If this 2-complex is planar, then the boundary of each cell adjacent to $\Id$ is a closed path from $\Id$ to itself that, read in the clockwise direction, is (a cyclic shift of) a relation in $R\cup\bar{R}$. We choose an arbitrary relation $g_1g_2\cdots g_k$ in $R$ and think of it as the first cell, drawn in the clockwise direction. So, we are fixing $g_1$ as the first generator in the cyclic ordering, while the second one (to the right in the plane and read \emph{from} the base vertex $\Id$) is $\bar{g_k}$. Now, to get the next generator to the right in the cyclic ordering, we have to choose a cyclic shift of a relation in $R\cup\bar{R}$ starting by $\bar{g_k}$, say $\bar{g_k}f_2\cdots f_l$. The next right generator will be $\bar{f_l}$. The condition is that $\bar{f_l}$ must not be a previously found generator.

\begin{figure}
\centering
\includegraphics[scale=0.7]{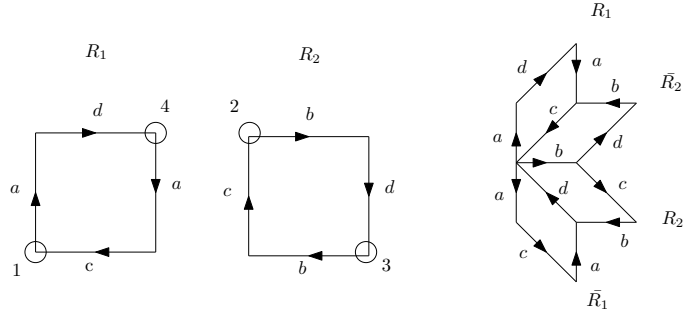}
\caption[fck]{A non geometric presentation. The circles numbered by $i$ indicate the angles used to attach the cell at step $i$ of the algorithm.}\label{co2}
\end{figure}

As an example, take for instance the presentation $P_1$ above. Set $R_1=adac$, $R_2=cbdb$. The process is illustrated in Figure~\ref{co2}. After three steps we get the generators $a,\bar{c},b,\bar{d}$, but in the fourth step it is not possible to continue: we cannot use the shift $\bar{d}\bar{b}\bar{c}\bar{b}$ of $\bar{R_2}$ (previously used in step 3) and we are forced to use the shift $\bar{d}\bar{a}\bar{c}\bar{a}$ of $\bar{R_1}$, what gives rise to $a$ as the next generator. Since $a$ was the starting generator and we have not yet completed the round (of all 4 generators and their inverses), it follows that the presentation $P_1$ cannot be geometric.

On the other hand, the previous procedure applied to the presentation $P_3$ above does give a complete cycle of all 4 generators and their inverses. The obtained cyclic ordering, together with all cells adjacent to the base vertex $\Id$, are shown in Figure~\ref{co3} (left).

\subsection{Computing the minimal bigons}
Once a cyclic ordering has been found for a presentation, the next step is to compute what we call the \emph{minimal bigons}
(see Proposition~\ref{Bigons}). Given two generators $x,y$ consecutive in the cyclic ordering, Proposition~\ref{Bigons}(b) states that there exists a unique bigon $\beta(x,y)$ of minimal length. Such a bigon is given by a pair of geodesics $\gamma_1,\gamma_2$ connecting the base vertex $\Id$ and another vertex $w$, with $\gamma_1$ starting at $\Id$ by the edge labelled $x$ and $\gamma_2$ starting at $\Id$ by $y$. Since both geodesics must have the same length, $\gamma_1\bar{\gamma_2}$ is a closed path of even length that, read clockwise, starts with the edge $x$ and ends with the edge $\bar{y}$. Recall that the \emph{length} of the bigon $\beta$ has been defined as the length of any of the two geodesics $\gamma_1,\gamma_2$ and is denoted by $k(x,y)$.

\begin{figure}
\centering
\includegraphics[scale=0.7]{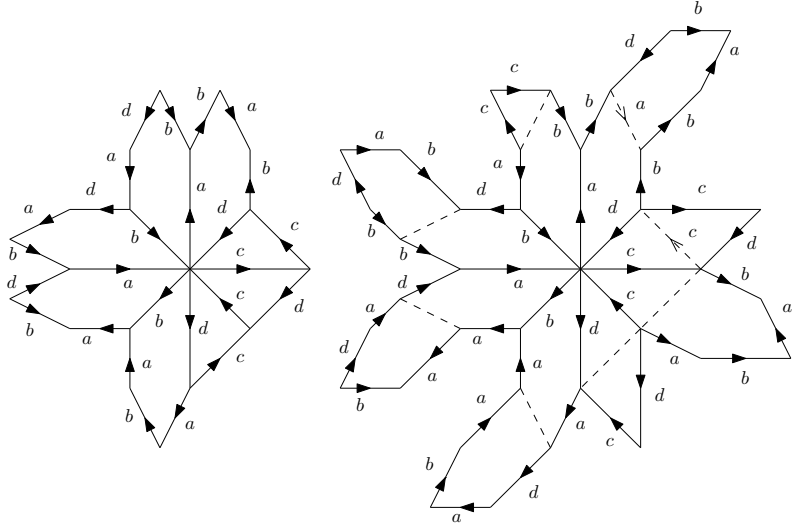}
\caption[fck]{Left: the cells adjacent to the base vertex $\Id$ and the cyclic ordering associated to the presentation $P_3$. Right: attaching cells to get the minimal bigons.}\label{co3}
\end{figure}

The existence of the minimal bigons is proved in Proposition~2.6 of \cite{L}. The proof is algorithmic and is implemented in the procedure {\bf MinimalBigons(R,co)}. In short and taking the above presentation $P_3$ as an example: consider the cell of the 2-complex $\textrm{Cay}^2(P_3)$ adjacent to $a$, whose boundary $B$ reads clockwise as $aba\bar{b}d$ (Figure~\ref{co3}, left). Since its length is odd, it cannot define a minimal bigon. So, we take its central symbol, $a$, look for the shift of a word in $R\cup\bar{R}$ starting by $a$, and glue together both cells along $a$. By Lemma~\ref{co}(c), there are exactly two of such shifts in $R\cup\bar{R}$, in this case $Y:=aba\bar{b}d$ and $Z:=a\bar{b}dab$. Since $Y$ is in fact a shift of $\bar{B}$, using $Y$ would yield a boundary containing the cancellation $b\bar{b}$, that cannot be geodesic. So, we are forced to use $Z$. From Lemma~\ref{co}(d) it easily follows that using $Z$ cannot produce cancellations. Once both cells are attached, in this case the obtained boundary has length 8, and thus defines a minimal bigon of length 4. See Figure~\ref{co3} (right). If the length of the boundary was odd, this procedure can be iterated and yields the minimal bigon after a finite number of steps.

\subsection{Choosing the cutting points and computing the kneading invariants}
The final stage of the whole process is the computation of the kneading matrix introduced in Section~\ref{milnor}. This task is assigned to the procedure {\bf KneadingMatrix(co,bigons,R)}. According to Theorem~\ref{mainTh}, the volume entropy of the presentation equals the topological entropy of the Bowen-Series-like map $\Phi_\Theta$ for \emph{any} choice of the set of cutting points $\Theta$ (note that specifying a particular map $\Phi_\Theta$ is equivalent to specifying $\Theta$). In consequence, we are free to choose the set of cutting points $\Theta$ in order for the kneading matrix of the corresponding map $\Phi_\Theta$ to be particularly simple. However, it is remarkable that a lot of information about any map $\Phi_\Theta$ in the family is independent of $\Theta$ and can be algorithmically extracted just from the combinatorial data encoded in the set of minimal bigons.

To illustrate it, we continue with the presentation $P_3$ above, whose minimal bigons are shown in Figure~\ref{co3} (right). The cyclic ordering reads clockwise as $\mathcal{O}=(a,\bar{d},c,\bar{c},d,b,\bar{a},\bar{b})$.

\begin{figure}
\centering
\includegraphics[scale=0.6]{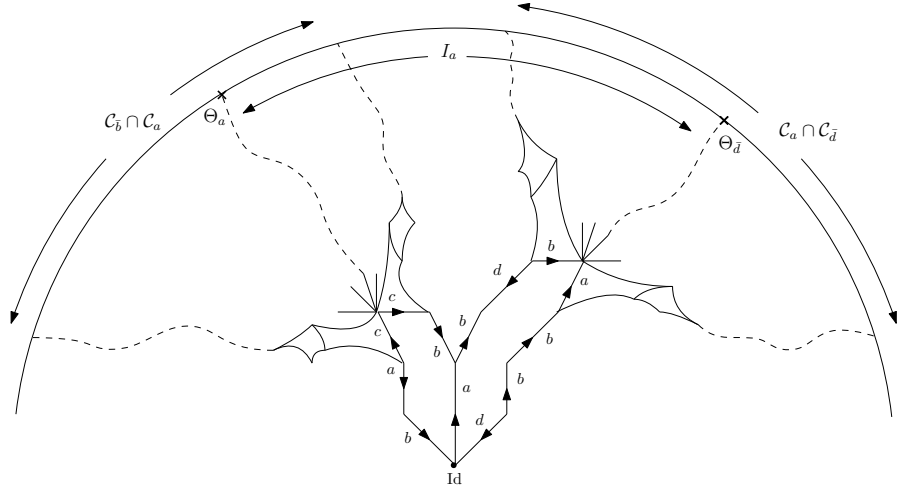}
\caption[fck]{The definition of the interval $I_a$ for the geometric presentation $P_3$.}\label{km1}
\end{figure}

Assume that for any pair $x,y$ of adjacent generators we have an (unspecified) cutting point $\theta_y$ in the intersection of the cylinders $\mathcal{C}_x\cap\mathcal{C}_y\subset S^1$, and set $I_x:=[\theta_x,\theta_y)\subset S^1$. See Figure~\ref{km1}. Let us see how to find the left and right images of any cutting point by a Bowen-Series-like map $\Phi$. As an example, the reader will find useful to get $\Phi(\theta_a)_{(-)}$ and $\Phi(\theta_a)_{(+)}$ following the picture shown in Figure~\ref{km1}. By the characterization of the intersection of cylinders given by Theorem~\ref{BamDia}, the cutting point $\theta_a$ has to be represented by two infinite words of the form $a\bar{b}\bar{c}\cdots$ and $\bar{b}\bar{a}c\cdots$. Considering $\theta_a$ as $a\bar{b}\bar{c}\cdots$ (a point of the cylinder $\mathcal{C}_a$) and taking into account that, by definition, $\Phi$ acts as a shift, it follows that $\Phi(\theta_a)_{(+)}$ is an infinite word starting by $\bar{b}\bar{c}\cdots$. In other words, $\Phi(\theta_a)_{(+)}\in I_{\bar{b}}$ by Lemma~\ref{ECcentral}(b). Analogously, $\Phi(\theta_a)_{(-)}$ can be read as an infinite word starting by $\bar{a}c\cdots$. Thus,
$\Phi(\theta_a)_{(-)}\in I_{\bar{a}}$.
%(for details see Lemma~\ref{ECcentralInterval})
As a quick recipe and after some practice, the reader will find easy to convince that, for any interval $I_x=[\theta_x,\theta_y)$, the right image of the left endpoint $\theta_x$ belongs to $I_z$, where $z$ is the generator following $x$ along the bigon that, read counterclockwise, starts by $x$, while the left image of the right endpoint $\theta_y$ belongs to $I_w$, where $w$ is the generator following $x$ along the bigon that, read clockwise, starts by $x$.

\begin{figure}
\centering
\includegraphics[scale=0.75]{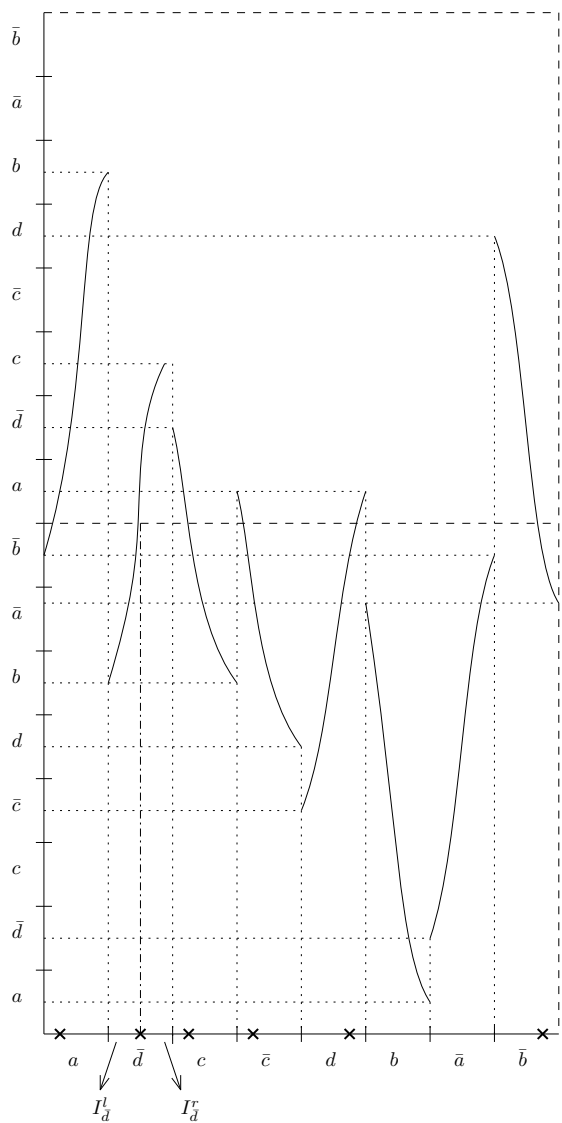}
\caption[fck]{A sketch of the graph of $\Phi$ for the geometric presentation $P_3$, where each interval $I_x$ has been labelled by $x$. The preimages of the leftmost point are marked with crosses in the horizontal axis. }\label{km2}
\end{figure}

At this moment we know which intervals contain the right and left images respectively of the left and right endpoints of any interval $I_x$. Since we know that $\Phi$ is a homeomorphism on each $I_x$, we can simply \emph{connect the dots} to get an approximate graph of $\Phi$. However, when connecting the dots, the question arises whether the map is increasing or decreasing on each interval. For instance, take the interval $I_c=[\theta_c,\theta_{\bar{c}})$. The previous recipe using the minimal bigons tells us that $\Phi(\theta_c)_{(+)}\in I_{\bar{d}}$ and $\Phi(\theta_{\bar{c}})_{(-)}\in I_b$. Now check Figure~\ref{co3}(right) to see that three clockwise consecutive edges $\bar{c},c,\bar{d}$ depart from the group element (vertex) $c$. Note that this three symbols are consecutive in the cyclic ordering $\mathcal{O}$, but \emph{in the counterclockwise direction}. In other words, the cyclic ordering of the edges departing from $c$ is the reverse of $\mathcal{O}$. It easily follows that the map has to be orientation reversing in $I_c$. See Figure~\ref{km2} for a sketch of a Bowen-Series-like map $\Phi$ corresponding to the presentation $P_3$, where the horizontal axis, that represents $S^1$, has been duplicated in the vertical direction. The reader can use again Figure~\ref{co3} to check that the edges departing from vertices $c,\bar{c},b,\bar{b}$ are organized clockwise as the reverse of the cyclic ordering $\mathcal{O}$, and so the map is decreasing on the corresponding intervals $I_c,I_{\bar{c}},I_b,I_{\bar{b}}$.

To apply the theory of Milnor-Thurston invariants, we have to consider the lifting $\tilde{\Phi}$ of $\Phi$ and the map $\map{\widehat{\Phi}}{[0,1]}$ defined by $\widehat{\Phi}(F)(x)=\tilde{\Phi}(x)-E(\tilde{\Phi}(x))$ where $E(y)$ denotes the integer part of $y$ (for simplicity and abuse of language, $\widehat{\Phi}$ will be denoted by $\Phi$ from now on). It is then clear that $\Phi$ is a piecewise monotone and discontinuous map. The set of points separating maximal intervals of monotonicity (\emph{turning points}) is precisely the set $\Theta$ of all cutting points \emph{plus} the set of preimages of the leftmost cutting point, $\theta_a$. There is at most one of such preimages, $m_x$, in any interval $I_x$, and in this case the interval $I_x$ splits into two subintervals separated by $m_x$, which we will call $I_x^l$ and $I_x^r$ (standing for \emph{left} and \emph{right}). See Figure~\ref{km2}. The ordered set of turning points for our presentation $P_3$ is then
\[ \theta_a<m_a<\theta_{\bar{d}}<m_{\bar{d}}<\theta_c<m_c<\theta_{\bar{c}}<m_{\bar{c}}<\theta_d<m_d<\theta_b<\theta_{\bar{a}}<\theta_{\bar{b}}<m_{\bar{b}}. \]

\begin{figure}
\centering
\includegraphics[scale=0.6]{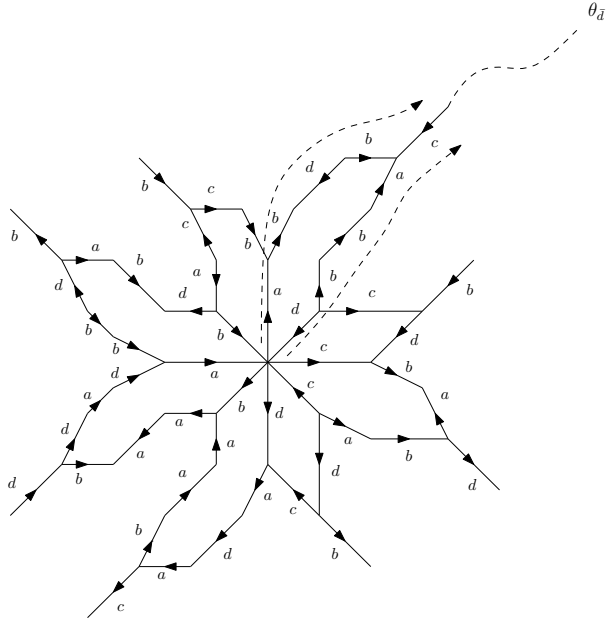}
\caption[fck]{Computing the left and right itineraries of the turning point $\theta_{\bar{d}}$.}\label{km3}
\end{figure}

Now we must choose the set of cutting points $\Theta$. Recall the notation $\{x_1,x_2,\ldots,x_{2N}\}$ introduced in Definition~\ref{standing} for the elements of the generating set. For each $1\le i\le 2N$, let $L_i^l$ and $L_i^r$ be the left and right geodesic segments of length $k(x_{i-1},x_i)$ defining the minimal bigon $\beta(x_{i-1},x_i)$. Let $\Theta$ be the \emph{middle parameter} according to the definition introduced in Section~\ref{middle}. This means that, for any $1\le i\le 2N$, $\theta_i$ is the $L_i^l$-middle point. To see how these notations fit to our example, consider for instance the generator $\bar{d}$. We have that
$(x_1,x_2,\ldots,x_8)=(a,\bar{d},c,\bar{c},d,b,\bar{a},\bar{b})$, so that $\bar{d}=x_2$. Then,
$\beta(x_1,x_2)=\beta(a,\bar{d})=\{L_2^l,L_2^r\}=\{ab\bar{d}b,\bar{d}bba\}$ (see Figure~\ref{km3}). Now we consider the centered continuation of $L_2^l=ab\bar{d}b$, which reads as $ab\bar{d}b\bar{c}\cdots$ since the opposite to $\bar{b}$ is $\bar{c}$. Although this is all the information we will need to proceed, one can check that in fact the centered continuation of
$L_2^l$ is $ab\bar{d}b(\bar{c}\bar{a}db)^\infty$. The cutting point $\theta_2=\theta_{\bar{d}}$ is the point of $S^1$ corresponding to this geodesic ray.

Now it's time to compute the kneading invariants for all turning points. To do it, we have to find the left and right itineraries of each turning point. We can
use Lemma~\ref{ECcentral}(a,b) to find the first $k$ symbols of the itinerary, where $k$ is the length of the corresponding bigon. For an example, consider again
the point $\theta_{\bar{d}}$. Moving along the minimal bigon $\beta(a,\bar{d})$, of length $k=4$ (see Figure~\ref{km3}), Lemma~\ref{ECcentral} states that the first four right iterates of $\theta_{\bar{d}}$ are
\[ \theta_{\bar{d}(+)}\in I_{\bar{d}},\ \Phi(\theta_{\bar{d}})_{(+)}\in I_b,\ \Phi^2(\theta_{\bar{d}})_{(+)}\in I_b,\ \Phi^3(\theta_{\bar{d}})_{(+)}\in I_a. \]
Analogously, the first four left iterates are
\[ \theta_{\bar{d}(-)}\in I_a,\ \Phi(\theta_{\bar{d}})_{(-)}\in I_b,\ \Phi^2(\theta_{\bar{d}})_{(-)}\in I_{\bar{d}},\ \Phi^3(\theta_{\bar{d}})_{(-)}\in I_b. \]

When an iterate belongs to an interval $I_x$ that splits into $I_x^l\cup I_x^r$, we can precise whether the iterate belongs either to $I_x^l$ or to $I_x^r$ simply taking into account the graph of the map shown in Figure~\ref{km2}. For instance, $\Phi^2(\theta_{\bar{d}})_{(-)}\in I_{\bar{d}}$ but, since $\Phi^3(\theta_{\bar{d}})_{(-)}\in I_b$, then $\Phi^2(\theta_{\bar{d}})_{(-)}$ belongs to the subinterval of $I_{\bar{d}}$ that is mapped onto $I_b$, which is $I_{\bar{d}}^l$. Analogously we have that
\[ \theta_{\bar{d}(-)}\in I_a^r,\ \Phi(\theta_{\bar{d}})_{(-)}\in I_b,\ \Phi^2(\theta_{\bar{d}})_{(-)}\in I_{\bar{d}}^l, \Phi^3(\theta_{\bar{d}})_{(-)}\in I_b, \]
\[ \theta_{\bar{d}(+)}\in I_{\bar{d}}^l,\ \Phi(\theta_{\bar{d}})_{(+)}\in I_b,\ \Phi^2(\theta_{\bar{d}})_{(+)}\in I_b,\ \Phi^3(\theta_{\bar{d}})_{(+)}\in I_a^l\cup I_a^r. \]

Now we should determine whether $\Phi^3(\theta_{\bar{d}})_{(+)}$ belongs either to $I_a^l$ or to $I_a^r$. This fact depends on the next symbol (the fifth one) in the infinite geodesic ray representing the cutting point $\theta_{\bar{d}}$. Since this symbol is $\bar{c}$, then $\Phi^3(\theta_{\bar{d}})_{(+)}\in I_a^r$.

See Figure~\ref{km3} for the complete set of symbols continuing the infinite rays at the vertices opposite to the identity in all minimal bigons. Now we can compute the first $k$ terms of the jump series. Going back again to the turning point $\theta_{\bar{d}}$, using the notation introduced in Section~\ref{milnor} and taking into account the increasing/decreasing character of $\Phi$ on each interval, we get
\[ \omega_0(\theta_{\bar{d}}^+)=I_{\bar{d}}^l,\ \omega_1(\theta_{\bar{d}}^+)=I_b,\ \omega_2(\theta_{\bar{d}}^+)=-I_b,\ \omega_3(\theta_{\bar{d}}^+)=I_a^r,\]
\[ \omega_0(\theta_{\bar{d}}^-)=I_a^r,\ \omega_1(\theta_{\bar{d}}^-)=I_b,\ \omega_2(\theta_{\bar{d}}^-)=-I_{\bar{d}}^l,\ \omega_3(\theta_{\bar{d}}^-)=-I_b.\]

Let us summarize what we have. Our algorithm has chosen a set of symbols $Z:=\{z_1,z_2,\ldots,z_{2N}\}$ such that $z_i$ is the opposite to the inverse of the last symbol in $L_i^l$. Note that $Z$ is totally determined by the minimal bigons. From this information and the previous paragraphs, we have been able to compute the first $k(x_{i-1},x_i)$ terms of the jump series $\Omega_{\theta_i}(t)$. By Lemma~\ref{essencial}, all the remaining terms vanish in all cutting points. Since, by Theorem~\ref{mainTh}, the volume entropy of the presentation is independent from $\Theta$, it follows that, from an algorithmic point of view, the information encoded in the minimal bigons is all we need to compute the volume entropy.

Here follows the complete list of kneading invariants for all turning points that results from our particular choice of the cutting points. For brevity, $I_x$, $I_x^l$ and $I_x^r$ will be respectively denoted by $x$, $x_l$ and $x_r$ from now on:

\begin{itemize}
\item[]
\item[] $\nu_{\theta_a}(t)=(a_l-\bar{b}_r)+(\bar{b}_l+\bar{a})t+(-\bar{c}_r+c_r)t^2$
\item[] $\nu_{m_a}(t)=(a_r-a_l)+t\nu_{\theta_a}(t)$
\item[] $\nu_{\theta_{\bar{d}}}(t)=(\bar{d}_l-a_r)+(-b+\bar{d}_l)t^2+(a_r+b)t^3$
\item[] $\nu_{m_{\bar{d}}}(t)=(\bar{d}_r-\bar{d}_l)+t\nu_{\theta_a}(t)$
\item[] $\nu_{\theta_c}(t)=(c_l-\bar{d}_r)+(-\bar{d}_l-c_r)t$
\item[] $\nu_{m_c}(t)=(c_r-c_l)+t\nu_{\theta_a}(t)$
\item[] $\nu_{\theta_{\bar{c}}}(t)=(\bar{c}_l-c_r)+(-a_r+b)t+(-b-\bar{a})t^2$
\item[] $\nu_{m_{\bar{c}}}(t)=(\bar{c}_r-\bar{c}_l)+t\nu_{\theta_a}(t)$
\item[] $\nu_{\theta_d}(t)=(d_l-\bar{c}_r)+(\bar{c}_r+d_l)t$
\item[] $\nu_{m_d}(t)=(d_r-d_l)+t\nu_{\theta_a}(t)$
\item[] $\nu_{\theta_b}(t)=(b-d_r)+(-\bar{a}-a_r)t+(-\bar{a}-d_r)t^2+(-\bar{b}_l-a_r)t^3$
\item[] $\nu_{\theta_{\bar{a}}}(t)=(\bar{a}-b)+(\bar{d}_l+a_r)t+(\bar{a}+a_l)t^2+(\bar{d}_r+\bar{b}_l)t^3$
\item[] $\nu_{\theta_{\bar{b}}}(t)=(\bar{b}_l-\bar{a})+(-d_l-\bar{b}_r)t+(-\bar{b}_r+\bar{b}_l)t^2+(\bar{a}-d_l)t^3$
\item[] $\nu_{m_{\bar{b}}}(t)=(\bar{b}_r-\bar{b}_l)+t\nu_{\theta_a}(t)$
\item[]
\end{itemize}

Finally, we formally write the above kneading invariants as a linear combination of the base
\[ (a_l,a_r,\bar{d}_l,\bar{d}_r,c_l,c_r,\bar{c}_l,\bar{c}_r,d_l,d_r,b,\bar{a},b_l,b_r) \]
and organize the coefficients of all invariants but the first one in matrix form, obtaining the following $13\times 14$ kneading matrix:

{\fontsize{9}{10}
\[ \arraycolsep=1.8pt \left( \begin{array}{cccccccccccccc}
-1+t & 1 & 0 & 0 & 0 & t^3 & 0 & -t^3 & 0 & 0 & 0 & t^2 & t^2 & -t \\
0 & -1+t^3 & 1+t^2 & 0 & 0 & 0 & 0 & 0 & 0 & 0 & -t^2+t^3 & 0 & 0 & 0 \\
t & 0 & -1 & 1 & 0 & t^3 & 0 & -t^3 & 0 & 0 & 0 & t^2 & t^2 & -t \\
0 & 0 & -t & -1 & 1 & -t & 0 & 0 & 0 & 0 & 0 & 0 & 0 & 0 \\
t & 0 & 0 & 0 & -1 & 1+t^3 & 0 & -t^3 & 0 & 0 & 0 & t^2 & t^2 & -t \\
0 & -t & 0 & 0 & 0 & -1 & 1 & 0 & 0 & 0 & t-t^2 & -t^2 & 0 & 0 \\
t & 0 & 0 & 0 & 0 & t^3 & -1 & 1-t^3 & 0 & 0 & 0 & t^2 & t^2 & -t \\
0 & 0 & 0 & 0 & 0 & 0 & 0 & -1+t & 1+t & 0 & 0 & 0 & 0 & 0 \\
t & 0 & 0 & 0 & 0 & t^3 & 0 & -t^3 & -1 & 1 & 0 & t^2 & t^2 & -t \\
0 & -t-t^3 & 0 & 0 & 0 & 0 & 0 & 0 & 0 & -1-t^2 & 1 & -t-t^2 & -t^3 & 0 \\
t^2 & t & t & t^3 & 0 & 0 & 0 & 0 & 0 & 0 & -1 & 1+t^2 & t^3 & 0 \\
0 & 0 & 0 & 0 & 0 & 0 & 0 & 0 & -t-t^3 & 0 & 0 & -1+t^3 & 1+t^2 & -t-t^2 \\
t & 0 & 0 & 0 & 0 & t^3 & 0 & -t^3 & 0 & 0 & 0 & t^2 & -1+t^2 & 1-t
\end{array} \right) \]}

Now we delete any column (for instance, the first one) and compute the determinant $D$ of the obtained $13\times 13$ matrix. The only factor of $D$ containing real roots in $[0,1)$ is
\[ t^{10}-3t^9-14t^8-13t^7-17t^6-12t^5-17t^4-13t^3-14t^2-3t+1, \]
and the smallest root is $\lambda\approx0.170554162$. From Theorem~\ref{kneading} it follows that the volume entropy of the presentation $P_3$ is $\log(1/\lambda)\approx \log(5.86324007)$.

We end this section providing the relevant computation details for the classical presentations for the orientable and nonorientable surfaces of genus 2, which are respectively
\[ P_4=\langle a,b,c,d\, |\, ab\bar{a}\bar{b}cd\bar{c}\bar{d}\rangle,\, \, P_5=\langle a,b,c,d\, |\, a^2b^2c^2d^2\rangle. \]

The above presentations are \emph{minimal} in the sense that the number of generators equals the rank of the group. Note first that the volume entropy of any minimal presentation has to be the same. Indeed,
the number of vertices at distance $m$ from the identity vertex depends only on the shape of the Cayley graph as defining a tiling of the plane. Since any minimal presentation of a surface of genus 2
has one single relation of length 8, the Cayley graph gives a tiling of the plane by octagons in all cases.

For $P_4$, the cyclic ordering is $(a,d,\bar{c},\bar{d},c,b,\bar{a},\bar{b})$. In this case all minimal bigons have length 4. The set of cutting points defines a partition
\[ \theta_a<m_a<\theta_{\bar{b}}<m_{\bar{b}}<\theta_{\bar{a}}<\theta_b<m_b<\theta_c<m_c<\theta_{\bar{d}}<\theta_{\bar{c}}<m_{\bar{c}}<\theta_d<m_d \]
and the map is orientation preserving in all intervals. The kneading matrix is

{\fontsize{9}{10}
\[ \arraycolsep=1.8pt \left( \begin{array}{cccccccccccccc}
-1+t & 1 & 0 & -t^3+t^4 & 0 & -t^4 & 0 & 0 & 0 & -t^2 & 0 & t^3 & t^2 & -t \\
0 & -1 & 1 & 0 & 0 & t^3 & t^2 & 0 & t & 0 & -t & -t^2 & 0 & -t^3 \\
t & 0 & -1 & 1-t^3+t^4 & 0 & -t^4 & 0 & 0 & 0 & -t^2 & 0 & t^3 & t^2 & -t \\
t^2 & -t^3 & 0 & -1 & 1 & -t & 0 & -t^2 & 0 & 0 & t & 0 & 0 & t^3 \\
t & 0 & 0 & -t^3+t^4 & -1 & 1-t^4 & 0 & 0 & 0 & -t^2 & 0 & t^3 & t^2 & -t \\
0 & -t^2+t^3 & 0 & 0 & t & -1 & 1 & -t & 0 & 0 & t^2-t^3 & 0 & 0 & 0 \\
t & 0 & 0 & -t^3+t^4 & 0 & -t^4 & -1 & 1 & 0 & -t^2 & 0 & t^3 & t^2 & -t \\
0 & -t & 0 & 0 & t^2 & 0 & t & -1 & 1 & 0 & -t^2+t^3 & -t^3 & 0 & 0 \\
t & 0 & 0 & -t^3+t^4 & 0 & -t^4 & 0 & 0 & -1 & 1-t^2 & 0 & t^3 & t^2 & -t \\
t & 0 & 0 & -t & 0 & -t^2 & -t^3 & 0 & 0 & -1 & 1 & t^3 & t^2 & 0 \\
0 & 0 & t & 0 & 0 & 0 & t^3 & 0 & t^2-t^3 & 0 & -1 & 1-t & 0 & -t^2 \\
0 & 0 & t^2-t^3 & 0 & 0 & 0 & 0 & 0 & t^3 & -t^2 & 0 & -1+t & 1 & -t \\
t & 0 & 0 & -t^3+t^4 & 0 & -t^4 & 0 & 0 & 0 & -t^2 & 0 & t^3 & -1+t^2 & 1-t
\end{array} \right) \]}

The relevant polynomial factor of the determinant obtained after deleting the first column is
\[ t^4-6t^3-6t^2-6t+1, \]
the well known growth polynomial for the minimal geometric presentations of surface groups \cite{AJLM,CW}, whose smallest root is about $0.143269846$. The volume entropy of presentation $P_4$ (and, thus, of $P_5$) is then about $\log(6.97983577)$. The same is true for the following minimal and geometric presentations $P_6$ and $P_7$ corresponding respectively to the orientable and non-orientable surfaces of genus 2. These presentations were called \emph{symmetric} in \cite{AJLM}:
\[ P_6=\langle a,b,c,d\, |\, abcd\bar{a}\bar{b}\bar{c}\bar{d}\rangle,\, \, P_7=\langle a,b,c,d\, |\, abcdcba\bar{d}\rangle. \]

It is worth noticing that previous computations \cite{AJLM,L} of the volume entropy for geometric presentations of surface groups have been obtained by using suitable Markov partitions for some particular Bowen-Series-like maps, and require to deal with determinants of much bigger matrices (for instance, a $57 \times 57$ matrix for the presentation $P_4$ above). We note also that the Cannon technique of \emph{cone types} introduced in \cite{Can84} and extended to all hyperbolic groups (see for instance \cite{GdlH}), allows in principle to compute the volume entropy but in general in a non algorithmic way.
It has been used in practice only for geometric presentations with all relations of constant even length $n$, giving rise to a regular tiling of the plane by $n$-gons. By contrast, the algorithmic approach given in this paper applies to any geometric presentation.

\bibliographystyle{plain}
\bibliography{biblio}
\end{document}